\documentclass[ts,nonblindrev]{informs3}

\OneAndAHalfSpacedXI
\usepackage{natbib}
 \bibpunct[, ]{(}{)}{,}{a}{}{,}%
 \def\bibfont{\footnotesize}%
\usepackage{bbm}
\def\N{\mathbb{N}}
\def\R{\mathbb{R}}

\def \E {\mathbb{E}}
\def \P {\mathbb{P}}

\def \-> {\rightarrow}

\def \argmin {\operatorname{argmin}}

\def \Var {\mbox{\bf Var}}
\def \cov {\mbox{\bf Cov}}

\usepackage{enumerate}
\usepackage{enumitem}
\usepackage{tabularx,booktabs,mathrsfs,setspace}%
\usepackage{mathtools}
\usepackage{multirow}
\usepackage{rotating}

\usepackage{float}
\usepackage{pgfplots}
\pgfplotsset{compat=1.9}
\usepackage{csvsimple}
\usepackage{makecell}
\usepackage{multirow}
\usepackage{adjustbox}
\usepackage[caption=false, position=top]{subfig}
\usepackage{alphalph}

\usepackage{comment}

\TheoremsNumberedThrough     % Preferred (Theorem 1, Lemma 1, Theorem 2)
\EquationsNumberedThrough    % Default: (1), (2), ...
\begin{document}
%%%%%%%%%%%%%%%%

%\setcounter{page}{1}
%\DOI{10.1287/xxxx.0000.0000}%
%\input{responseletter_EJOR_R1}

\RUNAUTHOR{Drent et al.}

\RUNTITLE{Dynamic Supply Mode Selection}

\TITLE{Efficient Emission Reduction Through Dynamic Supply Mode Selection}

\ARTICLEAUTHORS{
\AUTHOR{Melvin Drent}
\AFF{School of Industrial Engineering, Eindhoven University of Technology, Eindhoven, the Netherlands, PO BOX 513, 5600MB, \EMAIL{m.drent@tue.nl}}
\AUTHOR{Poulad Moradi, Joachim Arts}
\AFF{Luxembourg Centre for Logistics and Supply Chain Management, University of Luxembourg, Luxembourg City, Luxembourg, 6, rue Richard Coudenhove-Kalergi L-1359, \EMAIL{\{poulad.moradi, joachim.arts\}@uni.lu}}
} 

\ABSTRACT{%
Reducing the carbon footprint of global supply chains is a challenge for many companies. Governmental emission regulations are increasingly stringent, and consumers are increasingly environmentally conscious. Companies should therefore integrate carbon emissions in their supply chain decision making. 
In this paper, we study the inbound supply mode and inventory management decision making for a company that sells an assortment of products. Stochastic demand for each product arrives periodically and unmet demand is backlogged.
Each product has two distinct supply modes that
differ in terms of their carbon emissions, speed, and costs.
The company needs to decide when to ship how much using which supply mode such that total holding, backlog, and procurement costs are minimized while the emissions associated with different supply modes across the assortment remains below a target level.
We assume that shipment decisions for each product are governed by a dual-index policy for which we optimize the parameters. 
We formulate this decision problem as a mixed integer linear program that we solve through Dantzig-wolfe decomposition. We benchmark our decision model against two state-of-the-art approaches in a large test-bed based on real-life carbon emissions data. 
Relative to our decision model, the first benchmark lacks the flexibility to dynamically ship products with two supply modes while the second benchmark makes supply mode decisions for each product individually. Our computational experiment shows that our decision model can outperform the first and second benchmark by up to 15 and 40 percent, respectively, for moderate carbon emission reduction targets. 
}%

%\KEYWORDS{deterministic inventory theory; infinite linear programming duality; 
%  existence of optimal policies; semi-Markov decision process; cyclic schedule}

\KEYWORDS{Inventory; sustainability; policy comparison; multi-item; carbon emissions}
%\HISTORY{submitted january 2020}

\maketitle
%%%%%%%%%%%%%%%%%%%%%%%%%%%%%%%%%%%%%%%%%%%%%%%%%%%%%%%%%%%%%%%%%%%%%%

% Samples of sectioning (and labeling) in MSOM
% NOTE: (1) \section and \subsection do NOT end with a period
%       (2) \subsubsection and lower need end punctuation
%       (3) capitalization is as shown (title style).
%
%\section{Introduction.}\label{intro} %%1.
%\subsection{Duality and the Classical EOQ Problem.}\label{class-EOQ} %% 1.1.
%\subsection{Outline.}\label{outline1} %% 1.2.
%\subsubsection{Cyclic Schedules for the General Deterministic SMDP.}
%  \label{cyclic-schedules} %% 1.2.1
%\section{Problem Description.}\label{problemdescription} %% 2.

\section{Introduction}
\label{sec:intro}
The transportation sector has consistently been one of the most polluting European sectors for more than a decade now, and it is projected to remain so for the foreseeable future \citep{EEA}. 
%Projections for the transportation sector are particularly alarming. 
%Contrary to the energy sector, annual GHG emissions from the transportation sector are expected to increase over the next decade. %, even when taking into account numerous planned decarbonization initiatives. 
This, unfortunately, appears to be a trend that stretches beyond Europe. Recent analysis indicates that the G20 countries, currently responsible for 80\% of the global greenhouse gas (GHG) emissions, will see an increase of 60\% in their transportation sector emissions by 2050 \citep{ren21}. 
Prominent global climate targets, such as the ones outlined in the Paris agreement, will soon become unattainable \citep{EEA,emissionsGap2020}.

In light of the above, the European Union (EU) recently announced the Green Deal, a framework containing climate targets and policy initiatives that sets the EU on a path to reach carbon neutrality by 2050.
The Green Deal is legally enshrined in the European Climate Law, which states that member states are legally committed to meet the targets, and face penalties in case they do not meet these targets.
Being among the most polluting sectors, a key part of the Green Deal relates to policy initiatives that impact the transportation sector. For instance, the EU plans to extend the European emissions trading scheme (EU ETS) to include both road and maritime transport \citep{abnett_2020,EuropeanCommission3}.
Under the ETS, which until now includes only air transport, the EU enforces a cap on the total amount of GHG emissions from sectors covered by the scheme.
%Companies active in such a particular sector receive or acquire a certain amount of emission allowances that should cover its emissions over a given period; companies that emit more emissions than they have allowances for face heavy fines. 
The EU also investigates whether to increase fossil fuel taxation, thereby effectively raising the price of GHG emissions. 
More and more companies are reducing their emissions voluntarily as part of their corporate social responsibility. 
If not penalized by governments, companies that excessively pollute might still lose revenues as environmentally conscious customers take their business elsewhere \citep{dong2019}.

The developments described above highlight the urgency for companies to explicitly incorporate GHG emissions in their supply chain decision making. 
In this paper, we study the inbound supply mode and inventory decision making of a company that sells an assortment of products which are sourced from outside suppliers. The company wishes to keep the total GHG emissions associated with using different supply modes across the assortment below a certain target level in the most economically viable manner. 
As is often the case in practice, the company may rely on a third party logistics (3PL) provider for the inbound transport of the products. 3PLs typically offer several transport modes for the transportation of products -- these may differ in terms of transportation costs, transit times, and GHG emissions. Alternatively, there may be different suppliers (e.g. a near and offshore supplier) for a product that naturally have different costs and transportation emissions. In the remainder of the paper, we will use the terminology of a 3PL provider that offers multiple transportation modes.

The company can utilize the heterogeneity in the fleet of the 3PL to its advantage. While some transport modes are low emitting but slow, others may be fast but result in more emissions. Fast transport modes also typically come at the expense of a cost premium, and yet they are often relied upon when responsiveness is required (e.g. in case of imminent stock outs). 
Thus the company should rely dynamically on both transport modes. 
Implementing this holistically across the entire assortment of products allows the company to reduce emissions significantly for products for which it is relatively cost-efficient to do so and less for products for which this is more expensive. 
It additionally enables the company to reduce the total inventory and transportation costs by shipping the majority of products with a relatively cheap but slow transport modes while simultaneously resorting to faster but more expensive and often more polluting transport modes whenever expedited shipments are needed. 
While the advantages of dynamically selecting different transport modes are evident, two important and interrelated questions remain: 
\begin{enumerate}
\item When should the company ship how many units of which product with what transport mode?
\item What is the value of dynamically shipping products with different transport modes? 
\end{enumerate}
These questions are interesting but also intricate when one wishes to answer them for an entire assortment of products where the combined total of GHG emissions from transportation must not exceed a certain target level. 

To tractably answer the questions above, we focus on the setting where the 3PL offers two distinct transport modes for the transport of each product (or, equivalently, the setting where the company has already decided on the two transport modes for each product). 
These transport modes need not be the same for every product; they will depend on the characteristics of the suppliers as well as the 3PL (e.g., some products can be transported using aircraft or rolling stock while others can be shipped via inland waterways or ocean shipping).
We consider long distance and/or high volume transport lanes where the impact of the transport mode decisions of any particular individual shipper on the actual shipping and carbon footprint is negligible.
The company decides periodically how many units it wishes to transport with what transport mode and incurs mode specific unit transportation costs. 
Shipments arrive at the company after a deterministic transit time that depends on the transport mode that is used. Demand for each product in every period is stochastic and independent and identically distributed across periods. Any demand in excess of on-hand inventory is backlogged and satisfied in later periods. The company incurs per unit holding and backorder penalty costs, and the specific cost parameters may vary from product to product. The company seeks to minimize the long-run average holding, backorder, and transportation costs  while keeping the total long-run average GHG emissions from transportation of the entire assortment below a certain target level.

It is well-known that the optimal policy for the inventory system described above is complex, even in the simplest case of a single product and absent of the emission constraint \citep{whittemore1977optimal,feng2006base}. 
For the control of each product, we therefore use a heuristic policy that is originally due to \cite{veeraraghavan2008now}. They show that their so-called dual-index policy performs quite well compared to the optimal policy. The dual-index policy tracks two inventory positions for each product: The slow inventory position, which equals the on-hand inventory plus all in-transit products minus backlog, and the fast inventory position, which is defined similarly but includes only those in-transit products that are due to arrive within the transit time of the fastest transport mode. Under the dual-index policy, we place orders with both modes such that these inventory positions are kept at (or above) certain target levels, also referred to as base-stock levels. 
As such, the dual-index policy dynamically prescribes shipment quantities for both transport modes based on the net inventory level and the number of products that are still in-transit.
To find the optimal base-stock levels for the entire assortment of products, we formulate the decision problem as a non-linear non-convex integer programming problem. A partition reformulation of this problem allows us to use column generation techniques to solve the decision problem. These techniques enable us to decompose the complex multi-product decision problem into simpler sub-problems per product. Leveraging a separability result of \cite{veeraraghavan2008now}, we show that this sub-problem constitutes a special Newsvendor problem that can be solved efficiently through a simulation-based optimization procedure.
%We remark that while we assume a dual-index policy for each product, our model and analysis readily extends to a general class of policies where the fastest transport mode is operated according to the base-stock rule and the regular mode is operated according to any other rule.  

The main contributions of this paper are:
\begin{enumerate}
\item We are the first to study dynamic mode selection for an assortment of products with stochastic demand where the total average GHG emissions from the inbound transport of those products must be kept below a certain target level. 
\item We provide a tractable optimization model that finds a tight lower bound on the optimal solution as well as near-optimal feasible solutions within reasonable time. We show that our mathematical formulation of the decision problem allows us to decompose the non-linear non-convex integer programming problem into sub-problems per product. We leverage results from \cite{veeraraghavan2008now} to show that the sub-problems can be solved efficiently through a one dimensional search procedure in which each instance constitutes a Newsvendor type problem that is readily solved through simulation.
\item We perform an extensive computational experiment based on data from different industries. Through these experiments:
\begin{enumerate}[label=\roman*.]
    \item We establish the value of dynamic mode selection by comparing our model with a model in which only one transport mode per assortment product can be used. This value can go up to 15 percent in cost savings;
    \item We show that decomposing an aggregate carbon emission reduction target into targets for each product in the assortment individually is financially detrimental. Our holistic approach can lead to cost savings of over 40 percent relative to the approach with reduction targets per individual product;
    \item We find that the emission differences between transport modes relative to the cost difference between modes is the main determinant of emission reduction potential for a given assortment. In our experiments we find that, 20 percent of the products for which this ratio is highest contribute between 59 and 94 percent of the emission reduction.
\end{enumerate}
\end{enumerate}

%Transport mode selection, tactical, but we look at operational, two modes dynamically. 
%Static mode selection studied for single product with stochastic demand by \cite{Hoen} for an assortment product under determinisitic demand by \cite{Hoen2014Switching}, They show that substantial reductions in emissions can be achieved. However, static mode selection, i.e. one distinct mode for each product, reduces the flexibility to utilize other modes that may be preferable, particularly in the face of stochastic demand. For instance, the company may want to use a fast, but polluting mode in case of imminent backorders. Thus there is likely much value in dynamic mode selection over static mode selection. We establish this.  

The remainder of this paper is organized as follows. In Section \ref{sec:literatureReview}, we review the existing literature and position our work within the literature. Section \ref{sec:modelDM} contains the model description as well as the mathematical formulation of the decision problem. A column generation procedure to solve the decision problem is provided in Section \ref{sec:analysisDM}. We subsequently report on an extensive computational experiment in Section \ref{sec:numericalWorkDM}, and we provide concluding remarks in Section \ref{sec:conclusionDM}.

\section{Related literature}
\label{sec:literatureReview}
This paper integrates carbon emissions from inbound transportation into an inventory system with two supply modes. As such, our work contributes to the large stream of literature that studies  multi-mode or multi-supplier inventory systems. 
For an excellent overview of such systems, we refer the reader to the review papers of \cite{THOMAS2006245}, \cite{ENGEBRETHSEN20191} and \cite{svoboda2019review}, and references therein. We also contribute to the extensive body of literature that revolves around the integration of environmental aspects into supply chain decision making; see \cite{DEKKER2012671}, \cite{BRANDENBURG2014299}, \cite{BARBOSAPOVOA2018399}, and references therein, for an overview of this field.
In what follows, we focus on contributions that are most relevant to the present paper.

The decision how many products to order from which supplier is considered a canonical problem in the inventory management literature. It has been studied extensively since the sixties, mostly under the assumption that lead times are deterministic, that unmet demand is backlogged, and that only two distinct suppliers are at the disposal of the decision maker; the fastest being more expensive than the slowest.
\cite{fukuda1964optimal} and \cite{whittemore1977optimal} were the first to study this system. Assuming periodic review, they show that its optimal policy is a simple base-stock rule only under the assumption that the difference between the lead times of both suppliers is one period. For general lead time differences, the optimal policy is complex and can only be computed through dynamic programming for small instances. 
Since then, most researchers have focused on developing well-performing heuristic polices for which the best control policy parameters can be tractably obtained. %While there exists a large number of such heuristic
%policies, our focus here is on those policies most relevant to our work, and we refer the interested reader to \cite{svoboda2019review} for an excellent and recent overview of the entire field of multi-supplier inventory/distribution systems.

In this paper, we rely on the so-called dual-index policy to decide upon the shipment sizes for both transport modes for each product. Under this policy, which is originally due to \cite{veeraraghavan2008now}, two different inventory positions are kept track off: One that includes all outstanding shipments and one that includes only those outstanding shipments that are due to arrive within the lead time time of the fastest mode. 
\cite{veeraraghavan2008now} show numerically that the dual-index policy performs well compared to the optimal policy. In fact, \cite{MelvinArts} show that the dual-index policy is asymptotically optimal as the cost of the fastest transport mode and the backorder penalty cost become large simultaneously. The policy has received quite some attention in recent years  \citep[see, e.g.,][]{sheopuri2010new,arts2011efficient,sun2019robust}. 
We employ the dual-index policy because it is intuitive, has good performance, and can be optimized efficiently. 
Unlike the present work, the dual-index policy has so far been studied exclusively in single product settings under conventional cost criteria absent of any emission considerations. 

Within the transportation literature, inventory systems with multiple transport modes have received considerable attention too. 
To properly embed the present work in the existing literature, we group contributions to this field into two categories depending on the modelling assumptions regarding the usage of the available transport modes \citep[c.f.][]{ENGEBRETHSEN20191}. 
The first category, which we refer to as \textit{dynamic} mode selection, is concerned with inventory systems in which multiple transport modes are used simultaneously over a given (possibly infinite) planning horizon. 
Since we study an infinite horizon periodic review inventory model in which products can be transported with two distinct modes in each period, our work falls into this category -- as do all the inventory papers with two suppliers described so far. Only few papers exist in this category that explicitly account for carbon emissions, and the few that do differ substantially from the present work in terms of modelling choices as well as analysis. 
They either assume deterministic demand and a finite horizon \citep{PALAK2014198} or study the closely related yet different problem of splitting an order among several transport modes \citep{KONUR201746}. While not explicitly modeling carbon emissions, \cite{DONG201889} and \cite{lemmens2019} also study the benefit of dynamically switching between multiple transport modes in the context of multi-modal transport. They show that this can lead to more usage of less polluting transport modes without compromising on costs or responsiveness. 
Different from our work, all papers mentioned above consider the inventory control and transport mode decisions for a single product only. 

The second category concerns inventory systems in which a single transport mode is selected a priori at the start of a planning horizon; all replenishment orders until the end of that planning horizon are then shipped with this mode.  
We refer to this category as \textit{static} mode selection.
Two papers belonging to this category are particularly relevant to our work. 
\cite{Hoen2014Switching} study a periodic review inventory system under backlogging where inbound transport is outsourced to a 3PL that offers multiple transport modes. Assuming base-stock control for each mode, they are interested in selecting the transportation mode that leads to the lowest long-run average total cost consisting of holding, backlogging, ordering, and emission costs. For calculating transportation emissions, they rely on the well-known NTM methodology (we discuss this methodology in more detail later in Section \ref{subsec:descriptionDM} and Appendix \ref{sec:carbon}). 
We extend \cite{Hoen2014Switching} in two important directions. First, we move from static to dynamic mode selection, thereby incorporating the flexibility to dynamically switch between different transport modes for each product. Second, we consider an assortment of products under a single constraint on the total average transportation emissions from those products.
 \cite{Hoen} consider a similar constraint in a multi-product variant of the setting of \cite{Hoen2014Switching} under the assumption that demand is deterministic and inversely related to the price set by the decision maker. They show that because of the portfolio effect of such an assortment-wide emission constraint, carbon emissions from transportation can be reduced substantially at hardly any additional cost. 

The dual-index policy studied in this paper has the appealing feature that it can mimic static mode selection.
This is useful in our computational experiment where we establish the added value of dynamic mode selection over static mode selection. 
A closely related paper in that respect is \cite{berling2016} who study dynamic speed optimization of a single transport mode in a single-product stochastic inventory problem. 
They show that the value of dynamically controlling the speed of outstanding shipments, as opposed to a static speed policy, can be significant, both financially and from a carbon emission perspective.

Our review so far has almost exclusively revolved around papers on multi-period inventory systems. We note that there is also a stream of literature that integrates carbon emissions into single period multi-supplier models, see e.g., \cite{ROSIC2013109}, \cite{ARIKAN201415}, and \cite{CHEN2016196}. Similar to the majority of the papers discussed so far, these papers focus on single-product settings.

\section{Model description}
\label{sec:modelDM}
In this section, we first provide a description of
the inventory system under consideration and introduce the notation that we use throughout this
paper. We then describe the policy we propose to dynamically ship products with two transport modes. We conclude with providing a mathematical formulation of the decision problem. 

\subsection{Description and notation}
\label{subsec:descriptionDM}
We consider a company that sells an assortment of products. The inventories for these products are replenished from external suppliers through a third party logistics provider (3PL). A 3PL often offers several transport modes. We focus on the setting where the company has already decided upon two distinct transport modes that it would like to use for the transport of each product.
These two transport modes will differ in terms of costs, lead times, emissions, or a combination thereof. 
Given these two transport modes for each product, the operational question that remains is how many units of each product the company should transport using which transport mode at what time so that costs --holding, backlog, and ordering-- are minimized and an overall emission constraint is met. Companies will increasingly impose such constraints, either voluntarily or due to government regulation.  %citatie? Comment on multi-item vs single item?

The inventory system runs in discrete time with $t\in\N_0$ denoting the period index. Without loss of generality, we assume that the period is of unit length and coincides with the review epoch. Let $\mathnormal{J}=\{1,2,\ldots, |J| \}$ denote the nonempty set of products that the company offers for sale. Demand for product $j\in\mathnormal{J}$ across periods is a sequence of non-negative independent and identically distributed (i.i.d.) random variables $\{D_j^t\}$.
Any demand in excess of on-hand inventory is backlogged. Let $I^t_j$ denote the net inventory level (on-hand inventory minus backlog) of product $j$ at the beginning of period $t$ after any outstanding orders have arrived.
Each unit of product $j$ in on-hand inventory $(I^t_j-D^t_j)^+$ carried over to the next period and incurs a holding cost $h_j>0$.
Similarly, each unit of product $j$ in backlog $(D^t_j-I^t_j)^+$ incurs a penalty cost $p_j>0$. Here we use the standard notation $x^+ = \max(0, x)$. 

Each product can be shipped using two distinct transport modes from one supplier (or, equivalently, using one or two distinct transport modes from two distinct suppliers). Let $\mathnormal{M}=\{f,s\}$ denote the set of available transport modes, where we use $f$ and $s$ to refer to the faster and slower transport mode, respectively. Associated with the transport of one unit of product $j\in\mathnormal{J}$ with mode $m\in\mathnormal{M}$ is a cost $c_{j,m} \geq 0$, a deterministic lead time $l_{j,m} \in\N_0$, a distance traveled from the supplier to the company $d_{j,m}>0$, and a certain number of units CO2 emission $e_{j,m} \geq 0$. The weight of one unit of product $j$ is denoted $w_j>0$.
Recall that the company outsources its transport to a 3PL provider and hence has no control over the actual shipping.
We therefore consider variable emissions that depend only on product and transport mode specific characteristics as well as on distance traveled, and we refrain from incorporating a fixed emission factor per actual shipment.  This is a reasonable assumption for long distance and/or high volume transport lanes where the impact of the decisions of any particular individual shipper on the carbon footprint are negligible.
In line with previous literature that models transportation emissions in the context of mode selection \citep[e.g.,][]{Hoen,Hoen2014Switching}, %more citations here
we endow $e_{j,m}$ with the following structure which is based on the NTM methodology: %citation here
\begin{equation}\label{eq:emission}
    e_{j,m} = w_j(a_m+d_{j,m} b_m),
\end{equation}
where $a_m\geq 0$ and $b_m >0$ are a fixed and variable transport mode specific emission constant, respectively. 
Consistent with the NTM methodology, we assume that each product is shipped
with an averagely loaded transport mode. %Note that care should be taken with the units of the input parameters in \eqref{eq:emission} so that they are internally consistent. 
We define the lead-time difference between the fast and slow mode as $l_j = l_{j,s}-l_{j,f}\geq 0$ for each product $j\in\mathnormal{J}$.  
Conventional literature on dual model problems \citep[e.g.,][]{sheopuri2010new} imposes the assumption that the cost premium of using the fast mode does not exceed the lead-time difference multiplied by the penalty costs, i.e., $(c_{j,f}-c_{j,s})<l_{j}p_j$ to ensure that using the fast mode is attractive. We do not impose this assumption as whether using the fast supply mode will also depend on the target carbon reduction. When the fast supply mode is less polluting than the fast mode, the fast mode may become attractive even when  $(c_{j,f}-c_{j,s})<l_{j}p_j$. Conversely, the fast supply mode may become unattractive when the fast supply mode is more polluting, even when $(c_{j,f}-c_{j,s})<l_{j}p_j$.
%Observe that, contrary to the conventional literature on dual-mode problems \citep[e.g.,][]{sheopuri2010new}, we do not impose any assumptions on the unit transportation costs or emission units. 
Thus our model allows for situations where, e.g., the expensive transport mode is either the fastest and most polluting or the fastest and least polluting. 
Finally, the amount of items of product $j\in\mathnormal{J}$ to be shipped with transport mode $m\in\mathnormal{M}$ in period $t$ is denoted by $Q^t_{j,m}$. With this notation, observe that shipments $Q^{t-{l_{j,f}}}_{j,f}$ and $Q^{t-{l_{j,s}}}_{j,s}$ arrive in period $t$ so that we can write the following recursion for the inventory level $I_j^t$ of each product $j$:
\[
I_j^t = I_j^{t-1} - D_j^{t-1} + Q^{t-{l_{j,f}}}_{j,f} + Q^{t-{l_{j,s}}}_{j,s}.
\]

All notation introduced so far as well as notation that we will introduce later is summarized in Table \ref{tab:notation}.

\begin{table*}[!htbp]
  \caption{\textsf{Overview of notation.}}
	\fontsize{8pt}{9pt}\selectfont
\scriptsize
  \label{tab:notation}
  \begin{tabularx}{\textwidth}{l X}
    \toprule
    Notation& Description\\
		\midrule
	Sets & \\
    $\mathnormal{J}$ 			& Assortment; Set of all products.\\
    $\mathnormal{M}$ 			& Set of available transport modes, i.e. $\mathnormal{M}=\{f,s\}$. \\
	Input &  \\
		$D^t_{j}$				& Random demand for product $j\in\mathnormal{J}$ in period $t\in\N_0$.\\
		$p_j$ & Penalty cost for one unit of product $j\in\mathnormal{J}$ in backlog carried over to the next period. \\
		$h_j$ & Holding cost for one unit of product $j\in\mathnormal{J}$ in on-hand inventory carried over to the next period. \\
		$c_{j,m}$ & Cost of shipping one unit of product $j\in\mathnormal{J}$ with transport mode $m\in\mathnormal{M}$. \\
		$l_{j,m}$ & Transportation lead time for product $j\in\mathnormal{J}$ by transport mode $m\in\mathnormal{M}$.\\
	    $l_{j}$ &  Transportation lead time difference between the fast and slow mode for product $j\in\mathnormal{J}$, i.e. $l_j=l_{j,s}-l_{j,f}$.\\
		$d_{j,m}$ & Distance for the transport of product $j\in\mathnormal{J}$ with transport mode $m\in\mathnormal{M}$. \\
		$a_m$ & Fixed emission constant corresponding with transport mode $m\in\mathnormal{M}$.\\
		$b_m$ & Variable emission constant corresponding with transport mode $m\in\mathnormal{M}$.\\
		$w_j$ & Unit weight of product $j\in\mathnormal{J}.$ \\
				$e_{j,m}$ & Total units CO2 emission associated with shipping one unit of product $j\in\mathnormal{J}$ with transport mode
		$m\in\mathnormal{M}$.\\
		$\mathcal{E}^{max}$ & The maximally allowable carbon emissions for the transport of the entire assortment of products. \\ 
	Decision variables & \\
	    $S_{j,m}$ & Base-stock level for product $j\in\mathnormal{J}$ and transport mode $m\in\mathnormal{M}$. \\
	    $\Delta_{j}$ & Difference between the slow and fast base-stock level for product $j\in\mathnormal{J}$, i.e.  $S_{j,s}-S_{j,f}$. \\
		$\mathbf{S}_f$ & The vector $(S_{1,f},S_{2,f},\ldots,S_{\vert\mathnormal{J}\vert,f})$. \\
		$\mathbf{\Delta}$ & The vector $(\Delta_{1},\Delta_{2},\ldots,\Delta_{\vert\mathnormal{J}\vert})$. \\
	State variables & \\
	$I^t_{j}$ & Inventory level of product $j\in\mathnormal{J}$ at the beginning of period $t\in\N_0$ after orders have arrived. \\
	$IP^t_{j,f}$ & Fast inventory position of product $j\in\mathnormal{J}$  in period $t\in\N_0$ before shipping orders. \\
		$IP^t_{j,s}$ & Slow inventory position of product $j\in\mathnormal{J}$ in period $t\in\N_0$ after shipping orders with the fast transport mode. \\
		$Q^t_{j,m}$ & Amount of product $j\in\mathnormal{J}$ shipped with transport mode $m\in\mathnormal{M}$ in period $t\in\N_0$.\\
		$O^t_j$ & The overshoot of product $j\in\mathnormal{J}$ in period $t\in\N_0$, i.e. $(IP^t_{j,f}-S_{j,f})^+.$ \\
	Output of model &  \\ 
	$C(\mathbf{S}_f,\mathbf{\Delta})$ & Total long-run average holding, backlog, and ordering costs under a given control policy $(\mathbf{S}_f,\mathbf{\Delta})$.\\ 
	$\mathnormal{E}(\mathbf{S}_f,\mathbf{\Delta})$ & Total emissions from transportation under a given control policy $(\mathbf{S}_f,\mathbf{\Delta})$.\\  
	$C^{UB}_P (C^{LB}_P)$ & Upper (lower) bound for the optimal solution to Problem $(P)$. 
 \\
    \bottomrule
  \end{tabularx}
\end{table*}

\subsection{Control policy}
It is well-known that even for the simplest case where $|J|=1$ and absent of the emission constraint, the policy that prescribes the optimal shipment quantities is complex and can only be computed through dynamic programming for very small instances that are arguably not representative for practice \citep{whittemore1977optimal,feng2006base}. 
For the control of this inventory system, we therefore use a heuristic policy that is originally due to \cite{veeraraghavan2008now}. They show numerically that their so-called dual-index policy performs quite well compared to the optimal policy. The dual index policy tracks two indices: One that contains all orders that are still in-transit and one that contains only those in-transit orders that are due to arrive within the lead time of the fast mode. Based on these outstanding orders, the policy dynamically ships orders with both modes to keep these indices at certain target levels. In line with standard inventory management nomenclature, we also refer to these target levels as base-stock levels. More specifically, the policy operates as follows. At the beginning of every period $t$ after orders $Q^{t-{l_{j,f}}}_{j,f}$ and $Q^{t-{l_{j,s}}}_{j,s}$ have arrived, we review the \textit{fast} inventory position, which includes all in-transit orders -- i.e. shipped with both the slow and the fast transport mode -- that will arrive within the lead time of the fast transport mode: 
\[
IP^t_{j,f}=I^t_j + \sum_{k=t-l_{j,f}+1}^{t-1} Q^{k}_{j,f} + \sum_{k=t-l_{j,s}+1}^{t-l_j} Q^{k}_{j,s}.
\]
Then, if necessary, we place order $Q^t_{j,f}$ with the fast transport mode to raise the fast inventory position to its target level $S_{j,f}$. That is, the amount of product $j$ shipped in period $t$ with the fast transport mode equals:
\[ Q^t_{j,f} = (S_{j,f}- IP^t_{j,f})^+.\]
 After placing the fast shipment order, we inspect the \textit{slow} inventory position, which includes the fast order just placed
\[
IP^t_{j,s}=I^t_j + \sum\limits_{k=t-l_{j,f}+1}^{t} Q^{k}_{j,f} + \sum\limits_{k=t-l_{j,s}+1}^{t-1} Q^{k}_{j,s},
\]
and ship an order with the slow transport mode such that this inventory position is raised to its target level $S_{j,s}$, with $S_{j,s} \geq S_{j,f}$ since the fast inventory position is contained in the slow inventory position. Thus the amount of product $j$ shipped with the slow transport mode in period $t$ equals:
\[Q^t_{j,s} = S_{j,s}- IP^t_{j,s}.\]
Note that, contrary to the fast inventory position, the slow inventory position can never be larger than its base stock level $S_{j,s}$.
After shipping both orders, demand $D_j^t$ is satisfied or backlogged, depending on whether there is sufficient inventory available or not. 
The period then concludes with charging holding or backlog costs.

The order of events in a period $t$ for each product $j$ is thus as follows:
\begin{enumerate}
\item Orders $Q^{t-{l_{j,f}}}_{j,f}$ and $Q^{t-{l_{j,s}}}_{j,s}$ arrive with the fast and slow transport mode, respectively, and are added to the on-hand inventory $I^t_j$.
\item Review the fast inventory position and ship order $Q^t_{j,f}$ with the fast transport mode at unit cost $c_{j,f}$.
\item Review the slow inventory position and ship order $Q^t_{j,s}$ with the slow transport mode at unit cost $c_{j,s}$.
\item Demand $D^t_j$ occurs and is satisfied from on-hand inventory if possible, and otherwise backlogged.  
\item Incur a cost $h_j$ for any unit in on-hand inventory $(I_j^t-D^t_j)^+$ and a cost $p_j$ for any unit in backlog $(D^t_j-I^t_j)^+$.
\end{enumerate}
Observe that under a dual-index policy, slow orders entering the information horizon of the fast transport mode may cause the fast inventory position to exceed its target level. The amount by which the fast inventory position exceeds its target level is referred to as the overshoot. The fast inventory position in period $t$ after placing orders with both modes thus equals $S_{j,f} + O^t_{j}$, where $O^t_{j}$ denotes the overshoot for product $j$ in period $t$:
\[ O^t_j = IP^t_{j,f} + Q^t_{j,f} - S_{j,f} = (IP^t_{j,f} - S_{j,f})^+.
\]
Later, in Section \ref{subsec:sub-problem}, we shall see that computing the steady state distribution of the overshoot is crucial for determining the performance of a given control policy for a single product. 

We furthermore define $\Delta_j = S_{j,s} - S_{j,f}$, $j\in\mathnormal{J}$, so that the control policy for a product can be specified in terms of its base-stock levels  $S_{j,s}$ and $S_{j,f}$ or in terms of its base-stock level for the fast transport mode  $S_{j,f}$ and the difference  $\Delta_{j}$. 
We mostly use the latter specification in our subsequent analysis.
A control policy $(\mathbf{S}_f, \Delta)$ for the entire assortment of products consists of the vectors $\mathbf{S}_f=(S_{1,f}, S_{2,f},\ldots, S_{\vert J \vert, f})$ and $\mathbf{\Delta}=(\Delta_1, \Delta_2,\ldots, \Delta_{\vert J \vert})$. 

In what follows, for all sequences of random variables $X^t$, we define their stationary expectation as $\E[X] = \lim_{T\rightarrow\infty} (1/T) \sum\nolimits_{t=0}^T X^t$ and their distribution as $\P(X\leq x) =  \lim_{T\rightarrow\infty} (1/T) \sum\nolimits_{t=0}^T \mathbbm{1}\{{X^t\leq x}\}$, where $\mathbbm{1}\{ A \}$ is the indicator function for the event $A$. Whenever we drop the period index $t$ we refer to the generic stationary random variable $X$ with expectation and distribution as defined above. 
%We furthermore write $X(k)$ to denote the $k$-fold convolution of the random variable $X$. 
\subsection{Decision problem}
For a given control policy $(\mathbf{S}_f,\mathbf{\Delta})$, we define the total long-run average holding, backlog, and ordering costs per period for the entire assortment of products as
\begin{align} 
C(\mathbf{S}_f,\mathbf{\Delta}) & = \sum\limits_{j\in\mathnormal{J}} C_j(S_{j,f},\Delta_j) \nonumber \\  &=  \sum\limits_{j\in\mathnormal{J}}\left(h_j \E[(I_j-D_j)^+] + p_j \E[(D_j-I_j)^+] + \sum\limits_{m\in\mathnormal{M}} c_{j,m}\E[ Q_{j,m}]\right), \label{eq:cost}
\end{align}
and the total emissions as
\[\mathnormal{E}(\mathbf{S}_f,\mathbf{\Delta}) = \sum\limits_{j\in\mathnormal{J}} \mathnormal{E}_j(S_{j,f},\Delta_j) =  \sum\limits_{j\in\mathnormal{J}} \sum\limits_{m\in\mathnormal{M}} e_{j,m} \E[ Q_{j,m}], \]
where it is understood that the expectation operators are conditional on the control policy $(\mathbf{S}_f,\mathbf{\Delta})$. \cite{veeraraghavan2008now} show that $C(\mathbf{S}_f,\mathbf{\Delta})$ is well-defined for any control policy $(\mathbf{S}_f,\mathbf{\Delta})$ as long as $\E[D_j]<\infty$ for all products $j\in\mathnormal{J}$.

The objective of our decision problem is to minimize the total long-run average costs while keeping the total emissions below a target level $\mathcal{E}^{max}$.
Combining the above-mentioned leads to the following mathematical formulation of our decision problem which we refer to as problem $(P)$:
 \begin{alignat*}{5}
&(P)\qquad \qquad  	  &&\displaystyle \min       &&\quad && C(\mathbf{S}_f,\mathbf{\Delta}) &&  \\
&        				      &&\text{subject to} &&      \quad &&\mathit{E}(\mathbf{S}_f,\mathbf{\Delta}) \leq \mathcal{E}^{max}, &&\qquad \\
&											&&									&&	\quad		&&\mathbf{S}_f \in \mathbb{R}^{\vert J \vert},\quad \mathbf{\Delta}\in\mathbb{R}_0^{\vert J \vert}.  && 
 \end{alignat*}
Let ($\mathbf{S}_f^*,\mathbf{\Delta}^*$) denote an optimal solution to problem $(P)$ and let $C_P$ be the corresponding optimal cost.
Note that Problem $(P)$ is a non-linear non-convex knapsack problem where more than one copy of each item can be 
selected. It is well-known that even the simplest types of such knapsack problems are $\mathcal{NP}$-hard \citep[e.g.][]{Kellerer2004}. Since our knapsack is more complex, we conclude that Problem $(P)$ also falls in that same complexity class; it is hence likely that also for our problem no exact polynomial time solution algorithm exists. 

We remark that Problem $(P)$ enables companies to reduce carbon emissions from their inbound logistics by imposing a constraint on the maximally allowable carbon emissions. This is particularly useful for companies that seek to reduce carbon emissions proactively. However, companies may also take a reactive position and make supply mode decisions based only on inventory and transport costs. This cost will also include a carbon emission price component in regions where emissions are subject to carbon pricing mechanisms such as carbon crediting or taxing. We now briefly show that our model and analysis also apply to that setting. To that end, let $c^e$ denote the price of one unit of carbon emissions. This price can also depend on the transport mode $m$ and/or product $j$, but for ease of exposition we omit those dependencies. 
The long-run average costs per period in \ref{eq:cost} can now be redefined as follows:
\begin{align*} 
\tilde{C}(\mathbf{S}_f,\mathbf{\Delta}) &= \sum\limits_{j\in\mathnormal{J}}\left(h_j \E[(I_j-D_j)^+] + p_j \E[(D_j-I_j)^+] + \sum\limits_{m\in\mathnormal{M}} (c_{j,m}+c^ee_{j,m})\E[ Q_{j,m}]\right).
\end{align*}
The decision problem is now to minimize the long-run average costs per period:
 \begin{alignat*}{5}
&(\tilde{P})\qquad \qquad  	  &&\displaystyle \min       &&\quad && \tilde{C}(\mathbf{S}_f,\mathbf{\Delta}) &&  \\
&        				      &&\text{subject to} &&      \quad &&\mathbf{S}_f \in \mathbb{R}^{\vert J \vert},\quad \mathbf{\Delta}\in\mathbb{R}_0^{\vert J \vert}.  && 
 \end{alignat*}
 Problem $(\tilde{P})$ is less complex than the original decision problem since it does not involve a constraint that links the individual products. As such, Problem $(\tilde{P})$ can be decomposed in $|J|$ product specific problems, each of which can be solved individually. The column generation sub-problem that we will discuss in Section \ref{subsec:sub-problem} has a similar structure as the product specific problems of $(\tilde{P})$, and the solution method we discuss there thus also applies to $(\tilde{P})$.

\section{Analysis}
\label{sec:analysisDM}
This section focuses on finding the optimal control policy for Problem $(P)$. 
Our approach relies on the technique of column generation -- also named Dantzig-Wolfe decomposition after its pioneers \citep{dantzig1960decomposition}. 
This technique enables a natural decomposition of the original multi-product decision problem into smaller single-product problems that have more structure. 
We refer the interested reader to \cite{lubbecke2005selected} for a comprehensive survey on column generation.
Below, we first explain how we apply column generation to Problem $(P)$, and we then describe a simulation-based optimization method for solving the sub-problem of this column generation procedure.

\subsection{The column generation procedure}
\label{subsec:cgp}
We first reformulate decision problem $(P)$ as an integer linear program in which each binary decision variable corresponds to a certain combination of values for the decision variables of our original decision problem. We subsequently relax the integrality constraint and we call this problem the master problem $(\mathit{MP})$.
Formally, let $K_j$ be the set of all possible dual index policies for product $j\in\mathnormal{J}$. 
Each policy $k\in K_j$ is determined by its policy parameters $S^k_{j,f}$  and $\Delta^k_j$. Let $x^k_j\in\{0,1\}$ denote the decision variable that indicates whether policy $k\in K_j$ is selected ($x^k_j=1$) for product $j\in\mathnormal{J}$ or not ($x^k_j=0$). By relaxing the integrality constraint on this binary decision variable, we arrive at the mathematical formulation of the master problem $(\mathit{MP})$:
\begin{alignat}{5}
&(\mathit{MP}) \quad \quad && \displaystyle \min     && &&\sum\limits_{j\in\mathnormal{J}}\sum\limits_{k\in\mathnormal{K}_j}C_j(S^k_{j,f}, \Delta_j^k)  x_j^k &&  \label{problemMP}\\
&        			      &&\text{subject to} &&  &&\sum\limits_{j\in\mathnormal{J}}\sum\limits_{k\in\mathnormal{K}_j}\mathit{E}_j(S_{j,f}^k,\Delta^k_j)x_j^k \leq \mathcal{E}^{max},  &&  \label{eq:emissioncon}  \\
&											&&									&&			&& \sum\limits_{k\in\mathnormal{K}_j}x_j^k=1, && \quad \forall j\in\mathnormal{J} \label{eq:intcon} \\
&											&&									&&			&& \mbox{ }x_j^k \geq 0, && \quad \forall j\in\mathnormal{J},\forall k\in\mathnormal{K}_j. \nonumber
 \end{alignat}
Let $C^{\mathit{LB}}_P$ denote the optimal cost for master problem $\mathit{(MP)}$. Due to the linear relaxation of the integrality constraint on $x^k_j$, an optimal cost $C^{\mathit{LB}}_P$ is also a lower bound on the optimal cost for Problem $(P)$, $C_P$. 
 
Due to its large number of decision variables, master problem $(\mathit{MP})$ is solved using column generation. To this end, we first restrict master problem $(\mathit{MP})$ to a small subset $\tilde{K}_j \subseteq K_j$ of trivial policies per product $j\in J$ (i.e. columns) that are feasible for Problem $(P)$ (and thus also for Problem $(\mathit{MP})$). Such a trivial policy prescribes, for instance, to ship orders exclusively with the least polluting mode. 
This restricted problem is referred to as the restricted master problem $(\mathit{RMP})$. We then solve $(\mathit{RMP})$ to optimality, and we are interested in new policies $k\in K_j \setminus \tilde{K}_j$, $j\in\mathnormal{J}$, that will improve the objective value of $(\mathit{RMP})$ if they are added to $\tilde{K}_j$. 
Such policies $k\in K_j \setminus \tilde{K}_j$ are identified through solving a column generation sub-problem for each product $j\in J$. The objective function of such a sub-problem is the reduced cost as a function of the policy with respect to the current dual variables obtained through solving $(\mathit{RMP})$ to optimality. If a policy $k\in K_j \setminus \tilde{K}_j$ has a negative reduced cost, then adding that policy as a column to $\tilde{K}_j$ in $(\mathit{RMP})$ will reduce the objective value of $(\mathit{RMP})$.
More formally, the column generation sub-problem for product $j\in J$ has the following form:
  \begin{alignat*}{5} 
(\mathit{SP(j)})  ~  & ~ \displaystyle \min  ~    h_j \E[(I_j-D_j)^+] + p_j \E[(D_j-I_j)^+] + \sum_{m\in\mathnormal{M}} (c_{j,m} - \eta e_{j,m} )\E[ Q_{j,m}]  - \upsilon_j, \\
&        				      \text{subject to} \quad {S}_{j,f} \in \mathbb{R},\quad \Delta_j  \in \mathbb{R}_0,
 \end{alignat*}
 where $\eta$ denotes the dual variable of $(\mathit{RMP})$ that corresponds with emission constraint \eqref{eq:emissioncon} and  
 $\upsilon_j$ denotes the dual variable of $(\mathit{RMP})$ that corresponds with constraint \eqref{eq:intcon} that assures that for each product $j\in\mathnormal{J}$ a convex combination of policies is chosen. Note that these dual variables can also be interpreted as the Lagrange
multipliers of relaxing the corresponding constraints \citep{lubbecke2005selected}.
 If for product $j\in J$ a feasible solution for $(\mathit{SP(j)})$ exists with a negative objective value (i.e. a negative reduced cost), then this policy is added to $\tilde{K}_j$ since the objective value of $(\mathit{RMP})$ can be improved when solved with the enlarged set $\tilde{K}_j$. 
 
 We continue with iterating between optimizing $(\mathit{RMP})$ and finding new policies through solving $(\mathit{SP(j))}$, $j\in\mathnormal{J}$, until no product for which there is a policy with a negative reduced cost exists. 
 An optimal solution for $(\mathit{RMP})$ is then also an optimal solution for $(\mathit{MP})$. If this optimal solution contains integer values only, then it is also an optimal solution for $(P)$. If this is not the case, then we solve $(\mathit{RMP})$ one last time as an integer program to find an integer solution for $(\mathit{RMP})$, which is then also a feasible solution for $(P)$. Recent inventory literature has shown that solving the restricted master program as an integer program to arrive at an integer solution leads to good performance in terms of optimality gaps \citep[e.g.,][]{drent2020expediting, haubitz2020}, and often outperforms alternative approaches such as local searches or rounding procedures \citep{alvareze2013selective}. The corresponding cost of the resulting feasible solution is also an upper bound, denoted $C_{\mathit{UB}}^P$, for $C^P$.
 
 In the next section, we provide a simulation-based optimization procedure to solve the column generation sub-problem $(\mathit{SP}(j))$. %While this procedure solves $(\mathit{SUB}(j))$ to optimality, we note that any policy with a negative reduced cost can potentially improve the objective value of $(\mathit{RMP})$ if they are added; that is, not exclusively those policies that are both optimal for $(\mathit{SUB}(j))$ and have a negative reduced cost.

 \subsection{Solving the column generation sub-problem}
  \label{subsec:sub-problem}
 The column generation sub-problem $(\mathit{SP(j)})$ has the same structure as the problem studied by \cite{veeraraghavan2008now}. We follow their simulation-based optimization procedure to solve $(\mathit{SP(j)})$. 
 This procedure is grounded in the following separability result that allows us to find the optimal $S_{j,f}$ for a given $\Delta_j$ as the solution to a special Newsvendor problem.  
\begin{lemma}{\citep[Proposition 4.1]{veeraraghavan2008now}} 
\label{lem:seperability}
The distributions of the overshoot $O_j$, the fast transport mode shipment size $Q_{j,f}$, and the slow transport mode shipment size $Q_{j,s}$ are functions of $\Delta_j$ only, independent of $S_{j,f}$. 
\end{lemma}

 Let $O^t_j(\Delta_j)$ denote the overshoot of product $j$ in period $t$ for a given $\Delta_j$.
 Recall that the fast inventory position of product $j$ in period $t$ after shipping equals $S_{j,f} + O^t_j(\Delta_j)$. Consequently, for the net inventory level of product $j$ in each period $t$, we can also write
\begin{equation}
  I_j^{t}  = S_{j,f} - \left( \sum\nolimits_{k=t-l_{j,f}}^{t-1} D_j^k - O^{t-l_{j,f}}_j(\Delta_j) \right).  \label{eq:netInventory}
\end{equation}
By plugging \eqref{eq:netInventory} in the objective function of $(\mathit{SP(j)})$, and using Lemma \ref{lem:seperability} as well as the fact that in each period the overshoot is independent of the demand and $S_{j,f}$, we readily recognize that for given $\Delta_j$ the objective function is convex in $S_{j,f}$. This implies the following result. %rewrite as it's not just plugging in
\begin{lemma}\citep[Theorem 4.1]{veeraraghavan2008now}
\label{lem:optimalS}
The optimal base-stock level $S_{j,f}^{*}$ for a given $\Delta_j$, denoted $S_{j,f}^*(\Delta_j)$, equals
\[
S_{j,f}^{*}(\Delta_j)= \inf \left\{S_{j,f} \in \R : \P \left(   \sum\nolimits_{i=1}^{l_{j,f}+1} D_i - O_j(\Delta_j) \leq S_{j,f} \right) \geq \frac{p_j}{p_j+h_j} \right\}. 
\]
\end{lemma}

It now remains to calculate the objective value of $\mathit{SUB}(j)$ for given $\Delta_j$ and corresponding $S_{j,f}^*(\Delta_j)$. 
Observe that in each period immediately after shipping orders with both modes, the slow inventory position equals the fast inventory position plus the overshoot and all remaining outstanding slow orders. Since these inventory positions are equal to their respective base-stock levels following order placement, we have for each product $j$ in each period $t$:   
\begin{equation}\label{eq:totalOrders} S_{j,s} = S_{j,f}+O^t_j + \sum\limits_{k=0}^{l_j-1} Q^{t-k}_{j,s}.
\end{equation}
From \eqref{eq:totalOrders} it follows that $\E[Q_{j,s}] =     (\Delta_j - \E[O_j]) / l_j $. Since under backlogging the sum of both orders must on average be equal to the period demand, we finally find $\E[Q_{j,f}] = \E[D_j] - \E[Q_{j,s}]$.

To solve $(\mathit{SP(j)})$, $j\in\mathnormal{J}$, to optimality, it thus suffices to perform a one-dimensional search over $\Delta_j$. For each $\Delta_j$, we compute the stationary distribution of the overshoot. With this stationary distribution we readily find the optimal base-stock level  $S_{j,f}^*(\Delta_j)$ through Lemma \ref{lem:optimalS} and the total reduced cost through the identities following Equation $\eqref{eq:totalOrders}$. 
As there is in general no closed-form expression for the stationary distribution of the overshoot, we follow \cite{veeraraghavan2008now} and rely on simulation to compute this distribution. %(\cite{arts2011efficient} propose an approximate method to compute the stationary distribution of the overshoot using Markov chains, thus refraining from simulation.)

Note that our optimization model and analysis readily extends to  settings where the slow transport modes of all (or some) products are operated according to any other rule that depends only on the current overshoot as well as all in-transit orders that are not yet included in the fast inventory position. 
That is, any other rule that depends only on the information state $(O_j^t,Q^{t-1}_{j,s}, Q^{t-2}_{j,s},\ldots, Q^{t-l_{j}+1}_{j,s})$, $j\in\mathnormal{J}$, $t\in\mathbb{N}_0$. Most well-performing control policies satisfy this condition, e.g., the Capped Dual-Index policy \citep{sun2019robust}, the Tailored Base-Surge policy \citep{allon2010global}, and the Projected Expedited Inventory Position policy \citep{MelvinArts}.
\cite{sheopuri2010new} show that for such control policies, the stationary distribution of the overshoot is a function of only the parameter(s) for operating the slow transport mode, and that consequently a Newsvendor result similar to Lemma \ref{lem:optimalS} holds for all such policies. 

\section{Computational experiment}
\label{sec:numericalWorkDM}
Our test-bed has three different types of assortments of products, each representing a different type of industry. The first assortment type consists solely of products for which emissions from the
slowest transport mode are less than the emissions from the faster transport mode. This assortment is inspired by apparel goods that are delivered from Vietnam to Europe by sea transport as the slowest mode and by air transport as the fastest mode. In this example, the fast supply mode has a higher carbon footprint.
The opposite holds for the second assortment type. Here we are inspired by industrial goods that are delivered from China to Europe by sea transport as the slowest mode and from Germany by truck as the fastest mode. In this case, the slow supply mode is associated with more emissions as goods are transported over a longer distance. The third assortment has products of both types.

We perform a parametric computational experiment. The base case is set up as follows. 
We consider 100 products for each assortment type, i.e. $|J|=100$. 
The input parameters in the base case are identical for all three assortments, except for the carbon emissions from transportation.
For each product $j\in J$, the period demand $D_j$ follows a negative binomial distribution. To create heterogeneous assortments, the parameters of this negative binomial distribution are randomly drawn from two separate distributions for each product $j$. The mean $\mu_{D_j}:=\E[D_j]$ is randomly drawn from a gamma distribution with mean 100 and coefficient of variation of 0.5. 
The coefficient of variation $CV_{D_j}:=\sqrt{\Var[D_j]}/\mu_{D_j}$ is randomly drawn from a shifted beta distribution with mean 0.9, standard deviation of 0.25, and shifted to the right by 0.3. Since low demand products typically have higher holding cost, the holding cost $h_j$ is negatively correlated with the mean demand $\mu_{D_j}$ of each product $j$ through a Gaussian copula with a fitting covariance matrix. In particular, $h_j$ is drawn from a gamma distribution with mean 1 and coefficient of variation equal to $CV_{D_j}$, with a Pearson correlation coefficient of -0.5. Details regarding our approach to generate correlated random numbers are relegated to Appendix \ref{sec: random_numbers}.

We set $l_{j,f}$ and $c_{j,s}$ to 0 for all products, and focus on $l_{j,s}$ and $c_{j,f}$, which now coincide with the lead time difference and the cost premium of product $j\in J$, respectively. We set $l_{j,f}=3$ for all products.  
The back-order penalty cost $p_j$ for product $j$ is a function of its holding cost $h_j$. The ratio between the holding and the penalty cost is an important determinant of the service level in an inventory system. Therefore we set $p_j=\psi_j \chi_j^p h_j$ where $\psi_j$ is a parameter we use to control the ratio between $p_j$ and $h$, and $\chi_j^p$ is a random perturbation. That is, $\chi_j^p$ has a shifted beta distribution with mean 0.98, standard deviation 0.1, and shifted 0.02 to the right.
 
The cost premium $c_{j,f}$ of the fast transport mode of product $j$ equals $\chi_j^c p_j l_{j,f}$, where $\chi_j^c$ has a beta distribution with mean 0.25 and standard deviation 0.1. Table \ref{tab:baseCase} provides a summary of how we randomly generated the products of the base case; we shortly explain how we randomly generated the emission units for these products for all the three assortment types. 
In Table \ref{tab:baseCase}, $NB(\mu,cv)$ denotes a negative binomial random variable with mean $\mu$ and coefficient of variation $cv$, $\Gamma(\mu,cv)$ denotes a gamma random variable with mean $\mu$ and coefficient of variation $cv$, and $B(\mu,\sigma,s)$ denotes a beta random variable with mean $\mu$ and standard deviation $\sigma$ that is shifted to the right by $s$; if we drop $s$ then this beta random variable is not shifted, i.e. $B(\mu,\sigma,s) =_d B(\mu,\sigma)+s$ where $=_d$ denotes equality in distribution. 

\begin{table}[htbp]
  \centering
  \caption{\textsf{Generating the base case input parameters.}}
    \fontsize{8pt}{9pt}\selectfont
    \begin{tabular}{ll}
    \toprule
    \multicolumn{1}{l}{Input parameter} &       \multicolumn{1}{l}{Generation}  \\
    \midrule 
    $\vert \mathnormal{J} \vert$ & 100 \\
    $D_j$ &       $\mathnormal{NB}(\mu_{D_j},\text{CV}_{D_j})$, with $\mu_{D_j}\sim\Gamma(100,0.5)$ and   $\text{CV}_{D_j}\sim \mathnormal{B}(0.9,0.25,0.3)$ \\
    $l_{j,f}$     &        0 \\
    $l_{j,s}$     &        3 \\
    $h_j$     &         $\Gamma(1,0.5)$, with $\rho_{\mu_{D_j},h_j} = -0.5$ \\
    $p_j$     &         $\psi_{p} \chi_j^p h_j$, with $\psi_{p} = 9$ and $\chi_j^p\sim \mathnormal{B}(0.98,0.1,0.02)$ \\
    $c_{j,s}$     &        0 \\
    $c_{j,f}$     &       $\chi_j^c p_j l_{j,f}$, with $\chi_j^c \sim \mathnormal{B}(0.25,0.1)$ \\
    \bottomrule
    \end{tabular}%
  \label{tab:baseCase}%
\end{table}%

As explained in Section \ref{sec:modelDM}, we rely on the NTM framework \citep{NTMframework} to set emissions based on the structure of equation \eqref{eq:emission}. 
We apply the NTM framework to data from the UN Comtrade Database \citep{UNComtrade} to obtain sample emission units; details regarding this methodology are relegated to appendix \ref{sec:carbon}. 
We then apply maximum likelihood estimation on these sample unit emissions to obtain three distinct sets of two distribution functions; two for each assortment type of products. The number of carbon emission units $e_{j,m}$ for product $j\in J$ with transport mode $m\in M$ are then randomly drawn from these distributions. These distributions are presented in Table \ref{tab:emissionUnits}. In this table, $LN(\mu,\sigma)$ denotes a random variable whose logarithm is normally distributed with mean $\mu$ and standard deviation $\sigma$, and $WB(\lambda,k)$ denotes a Weibull random variable with scale $\lambda$ and shape $k$.  

\begin{table}[htbp]
  \centering
  \caption{\textsf{Generating the emission units for the base case.}}
    \fontsize{8pt}{9pt}\selectfont
    \begin{tabular}{lll}
    \toprule
    Assortment type & \multicolumn{1}{l}{Emission parameter} & \multicolumn{1}{l}{Generation} \\
    \midrule
    \multirow{2}[1]{*}{1} &   $e_{j,s}$    &  $\Gamma(0.35,0.21)$ \\
          &   $e_{j,f}$    &  $e_{j,s}$ + $\mathnormal{LN}(1.52,0.21)$\\
          \addlinespace
    \multirow{2}[0]{*}{2} &   $e_{j,f}$    &  $\Gamma(0.19,1.27)$ \\
          &   $e_{j,s}$    &  $e_{j,f}$ + $\Gamma(2.19,1.27)$\\
          \addlinespace
    \multirow{2}[0]{*}{3} &   $e_{j,f}$    &  $\mathnormal{WB}(0.87,0.77)$ \\
          &   $e_{j,s}$    &  $\Gamma(3.31,1.34)$\\
    \bottomrule
    \end{tabular}%
  \label{tab:emissionUnits}%
\end{table}%

The total allowable carbon emissions from transportation $\mathcal{E}^\text{max}$ is set as a percentage of the total \textit{reducible} carbon emissions. 
For each instance of the test-bed, the reducible carbon emissions is defined as the difference between the total amount of carbon emissions of the control policy that is optimal for Problem $(P)$ absent of the emission constraint and the total amount of carbon emissions of the control policy that leads to the lowest possible total carbon emissions. 
The latter implies that each product is only shipped with its least polluting transport mode. Under the dual-index policy, setting $\Delta_j$ to zero implies that all orders for product $j\in J$ are shipped with its fastest transport mode. Alternatively, letting $S_{j,f}$ go to $-\infty$ implies that all orders for product $j$ are shipped with its slowest transport mode \citep[][]{veeraraghavan2008now}.

To evaluate the effectiveness of the column generation procedure in solving Problem $(P)$, we compute for each instance of the test-bed the relative difference between the total average cost under a feasible solution and the corresponding lower bound. That is,
\[ \%\mathit{GAP} = 100\cdot \frac{C^\mathit{UB}_P- C^\mathit{LB}_P}{ C^\mathit{LB}_P}, \]
where $C^\mathit{LB}_P$ and $C^\mathit{UB}_P$ are obtained using the methods described in Section \ref{sec:analysisDM}. In what follows, we also refer to this feasible solution as dynamic mode selection (DMS). Hence the long run average cost of dynamic mode selection equals $C^\mathit{UB}_P$.  

To quantify the benefit of using two transport modes dynamically rather than relying statically on one transport mode, we define for each instance of the test-bed a benchmark instance in which we can only select one transport mode for each product. As described above, the dual-index policy can mimic static mode selection (SMS). 
Hence, to find a feasible static mode selection solution to this benchmark instance of Problem $(P)$, we apply our column generation procedure of Section \ref{sec:analysisDM} in which we restrict the solution space for each product $j\in\mathnormal{J}$ such that all orders are shipped with either the fastest or the slowest transport mode. The mathematical formulation for the static mode selection approach as well as the benchmark approach described in the next paragraph can be found in Appendix \ref{sec:BM_approaches}. (Note that in the Master Problem $(\mathit{MP})$ for this approach the set of possible policies $K_{j}$ for each product $i\in\mathnormal{J}$ contains only two single transport mode policies.) The long run average cost of this solution is denoted $C^{\mathit{SMS}}_P$. To quantify the value of dynamic mode selection, we compare for each instance of the test-bed the long run average cost of static mode selection with the long run average cost of dynamic mode selection. That is, 
\[\%\mathit{SMS} = 100\cdot \frac{C^\mathit{SMS}_P- C^{UB}_P}{ C^\mathit{UB}_P}, \]
where $\% SMS$ indicates the relative increase in the long run average cost when the company chooses to rely on only one transport mode for each product in meeting an assortment wide constraint on total emissions.  

To quantify the portfolio effect, we define for each instance of the test-bed an additional benchmark in which we enforce an emission constraint for each product $j\in J$. 
This emission constraint is set as a percentage of the total reducible emissions of the individual product rather than of the entire assortment of products. 
We compute a feasible solution for this benchmark instance with the column generation procedure of Section \ref{sec:analysisDM}, which can readily be modified so that it can be applied to settings where we have additional emission constraints. (An alternate solution procedure is to use a simple search procedure per product.)
This approach, which we refer to as blanket mode selection (BMS), seems a plausible approach for most practitioners. Indeed, they can consider each product individually absent of the complicating linking emission constraint and yet they are guaranteed that the total emissions of the entire assortment is kept below the target level. 
The long run average cost of the blanket mode selection approach is denoted $C^{BMS}_P$. To quantify the portfolio effect, we compare for each instance of the test-bed the long run average cost of blanket mode selection with the long run average cost of dynamic mode selection.
That is, 
\[\%\mathit{BMS} = 100\cdot \frac{C^{BMS}_P- C^{UB}_P}{ C^{UB}_P}, \]
where $\%\mathit{BMS}$ indicates the relative increase in the long run average cost when the company enforces an emission constraint on each individual product rather than one single constraint for the emissions from the entire assortment. 

In solving the column generation sub-problem for each product, we simulated 10 samples of 9500 time periods following a 5000 time periods warm-up . The width of the 95 percent confidence interval of the long run average cost per period for each product was no larger than 3 percent of its corresponding point estimate for each instance of the column generation sub problem that we solved. The average computational time of our column generation procedure was 23 minutes, 5 minutes for blanket dynamic mode selection, and less than 5 sec for static mode selection.

\subsection{Results for the base case}
\label{subsec:base}
\begin{figure}[htb!]
	\centering
  \includegraphics[width=0.85\textwidth]{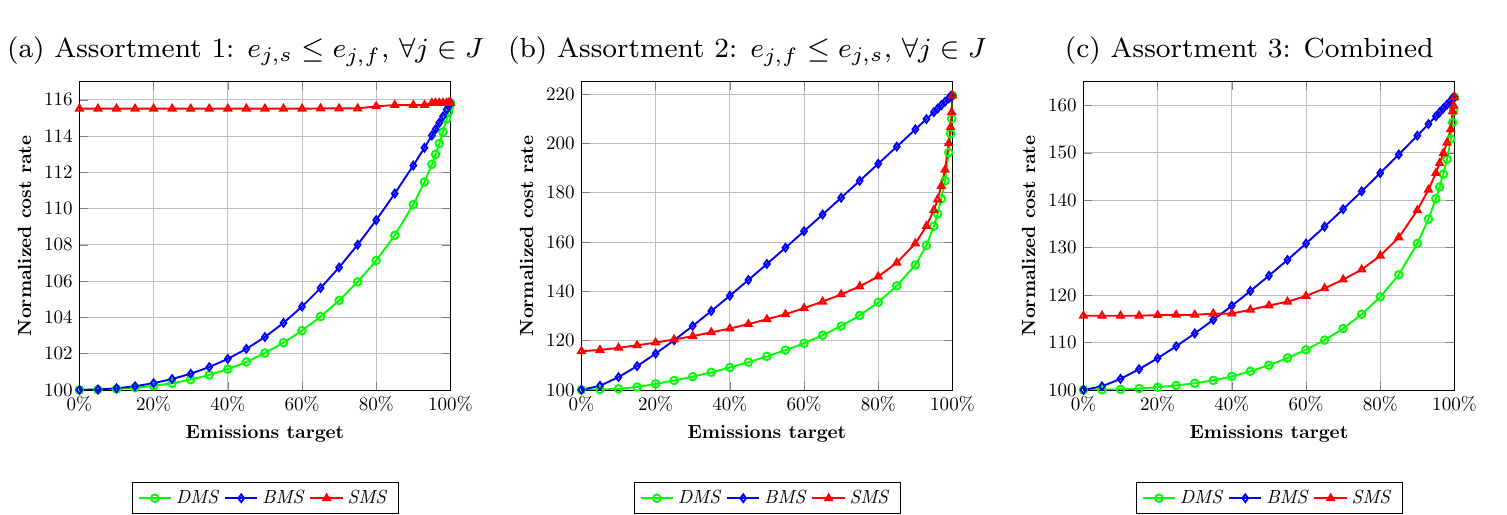}\\
	\caption{\textsf{Optimal normalized absolute costs of each approach for different targets on the reducible emissions.}}
  \label{fig:baseEfficient}
\end{figure}

Figure \ref{fig:baseEfficient} presents the normalized optimal average costs of each approach for each assortment group under emission targets that range from from 0 to 100 percent of the total reducible carbon emissions. Observe that for each assortment group, all approaches have the same performance when the emission target is set at 100 percent of the total reducible emissions. 
In this case, all approaches solely utilize the least polluting transport mode. Alternatively, when we impose no target on the emissions from transportation,
then both dual mode approaches perform equally well while the static mode selection approach seems to perform the poorest over all possible emission targets. Indeed, static mode selection is around 15 percent more expensive than both dual dual mode approaches for all assortment types when transportation emissions are not constrained. 
Based on Figure \ref{fig:baseEfficient}, we conclude that dynamic mode selection, as opposed to static and blanket mode selection, has great potential to efficiently curb carbon emissions from transportation at relatively little additional costs across all assortment types.  

%These observations are not surprising. 
%However, for moderate levels of emission targets, the results are less obvious and need further elaboration.  

\begin{figure}[htb!]
	\centering
  \includegraphics[width=0.85\textwidth]{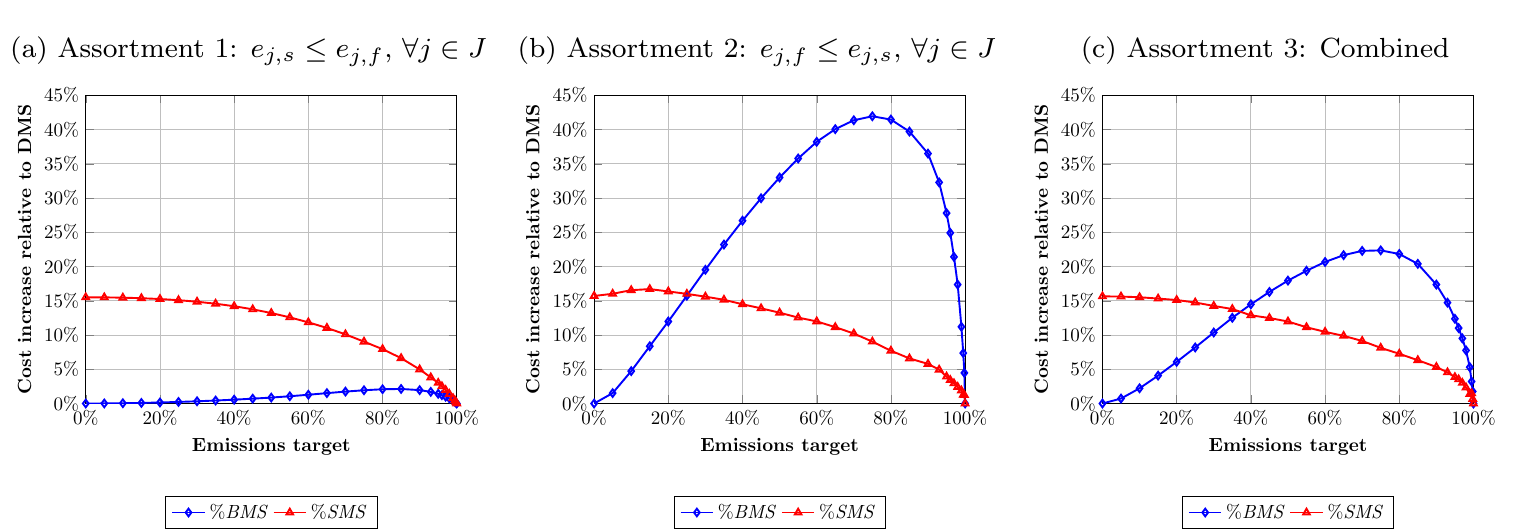}\\
	\caption{\textsf{Relative surplus of the optimal cost of the alternative approaches (BMS and SMS) compared to the DMS approach for different targets on the reducible emissions.}}
  \label{fig:base}
\end{figure}

We explicitly compare our dynamic mode selection with the benchmark approaches in Figure \ref{fig:base}, which presents the $\%\mathit{SMS}$ and $\%\mathit{BMS}$ percentages for each assortment group under emission targets that range from from 0 to 100 percent of the total reducible carbon emissions.
The figure indicates that the performance of static mode selection over dynamic mode selection is consistent across all assortments. 
The relative increase in its total cost over dynamic mode selection is the largest when there is no emission target, and gradually decreases as the emission constraint tightens. At moderate carbon emission targets, around 40 to 60 percent of the total reducible emissions, static mode selection still leads to increases in the total average cost per period of around 10 to 15 percent for all assortment types.    

The performance of the blanket mode selection approach depends on the specific assortment type. 
Figure \ref{fig:base}(a) illustrates that when the unit emissions from the fast transport mode are more than those from the slow transport mode, the performance of the blanket mode selection approach seems to be quite reasonable. 
This can be explained as in this setting, the cost of the fastest transport mode is larger than the cost of the slowest transport mode. 
The most polluting transport mode is thus also the most expensive transport mode. The portfolio effect is limited for this assortment type. 

For the other two assortment types, however, the cheapest transport mode is not necessarily also the least polluting transport mode, and the portfolio effect is more prevalent.
Figures \ref{fig:base}(b) and \ref{fig:base}(c)  show that $\%\mathit{BMS}$ can be more than 35 and 20 percent in assortment type 2 and 3, respectively, under carbon emission reduction targets of 50 percent. The static mode selection approach, which also takes advantage of the portfolio effect, even outperforms the blanket mode selection approach for quite some emission targets.  

\begin{figure}[htb!]
	\centering
    \includegraphics[scale=.95]{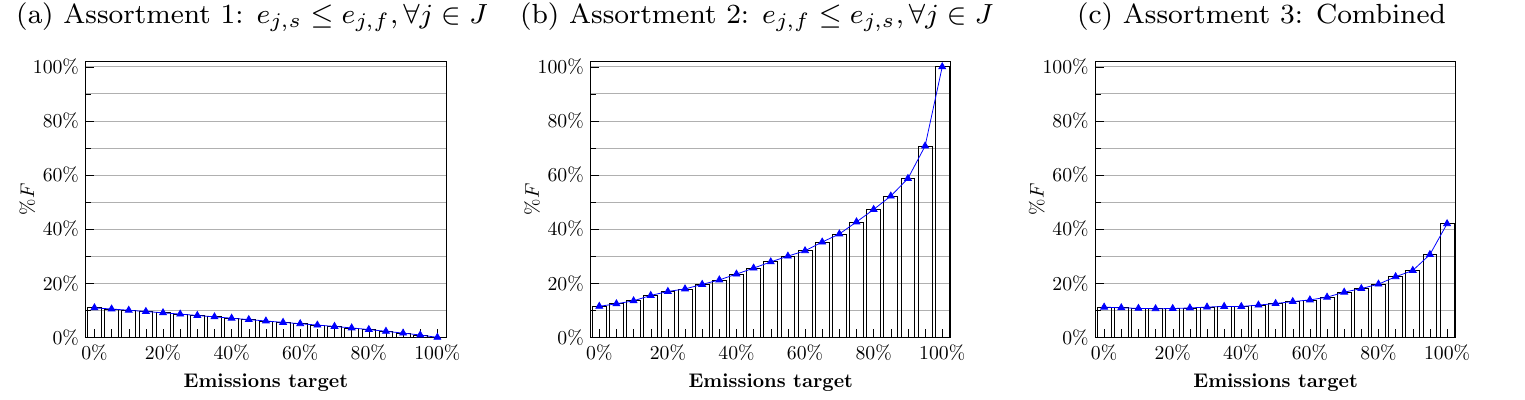}
    \caption{\textsf{Percentage of products shipped with the fast mode for different carbon emission reduction targets.}}
    \label{fig:Qsplitbar}
\end{figure}

Figure~\ref{fig:Qsplitbar} shows the usage of the fast supply mode as a function of the carbon emission reduction target under dynamic mode selection as measured by
\[
\%F = 100 \cdot \frac{1}{|J|}\sum_{j\in J} \frac{\E[Q_{j,f}]}{\E[Q_{j,f}] + \E[Q_{j,s}]}.
\]
Dynamic mode selection is economically attractive regardless of emissions targets. Carbon reduction targets make dynamic mode selection even more attractive for assortment type 2, and mixed assortments, but not for assortments of type 1.

% \begin{figure}[htb!]
% 	\centering
%     \begin{subfloat}[Assortment 1: $e_{j,s} \leq e_{j,f}, \forall j \in J$]
%         {
%         \centering
%         \includegraphics[scale=0.65]{PlotBarh1.pdf}
%         }
%     \end{subfloat}
%     \hspace{-6mm}
%         \begin{subfloat}[Assortment 2: $e_{j,f} \leq e_{j,s}, \forall j \in J$]
%         {
%         \centering
%         \includegraphics[scale=0.65]{PlotBarh2.pdf}
%         }
%     \end{subfloat}
%     \hspace{-5mm}
%     \begin{subfloat}[Assortment 3: Combined]
%         {
%         \centering
%         \includegraphics[scale=0.65]{PlotBarh3.pdf}
%         }
%     \end{subfloat}\\
%     \caption{\textsf{Percentage of products shipped with the fast mode for different carbon emission reduction targets.}}
%     \label{fig:Qsplitbar}
% \end{figure}

To recapitulate, the value of dynamically shipping products with two transport modes simultaneously rather than statically selecting one transport mode a priori is quite large. Regardless of the assortment type, $\%\mathit{SMS}$ is in between 5 and 15 percent for emission reduction targets up to 90 percent.
The portfolio effect depends on the specific assortment type. If the least polluting transport mode of each product is also its cheapest transport mode, then the fastest and most polluting transport modes are typically only relied upon in case of imminent backorders. This behavior remains in case of an assortment-wide emission target, and the portfolio effect is consequently rather limited.  
If the least polluting transport modes are not necessarily the cheapest transport modes, then there is substantial value to be reaped in optimizing the 
assortment of products under a single emission constraint rather than under separate emission constraints for each individual product. Indeed, $\%\mathit{BMS}$ can go up to 40 and 20 percent for assortment type 2 and 3, respectively.

Table \ref{tab:emissionsgap} below presents the average relative slack in the emission constraints for each assortment over the different emission targets considered in the base case analysis. The table shows that due to the binary nature of the static mode selection approach, the total average emissions under this approach are often substantially lower than the target level. This leads to particularly poor performance for assortment type 1. We observed in our computational experiments that for this assortment type, the static mode selection approach selects the cheapest and thus least polluting transport mode for almost all products under each emission target. 

\begin{table}[htbp]
  \centering 
	\caption{\textsf{Average slack in emission constraints for the base case analysis.}}
  \fontsize{8pt}{9pt}\selectfont
    \begin{tabular}{lrrr}
    \toprule
          & \multicolumn{3}{c}{Assortment type} \\
\cmidrule{2-4}    Approach & \multicolumn{1}{c}{1} & \multicolumn{1}{c}{2} & \multicolumn{1}{c}{3} \\
    \midrule
    DMS   & 0.00\%  & 0.00\%  & 0.00\% \\
    BMS   & 0.08\%   & 0.61\%  & 0.12\% \\
    SMS   & 12.16\% & 0.17\%  & 0.97\% \\
    \bottomrule
    \end{tabular}%
\label{tab:emissionsgap}
\end{table}%

%Seemingly, other better performing dual-sourcing heuristics such as the projected expedited inventory position heuristic \citep{MelvinArts} can reach the same targets at even lower costs.

The average $\%GAP$ of the base case over all emission targets is less than 0.01 percent, indicating that the column generation procedure finds feasible solutions that are close to optimal. Such a low average $\%GAP$ occurs because there can be at most 1 product for which the optimal solution to Problem $(MP)$ is fractional. Indeed, Problem $(MP)$ has $|J| + 1$ constraints and an optimal solution for this problem has the same number of basic variables.
Constraint \eqref{eq:intcon} assures that for each product $j\in J$ a convex combination of policies is chosen. As such, there is at least one basic variable for each product $j$. This implies that there is at most 1 product for which the optimal solution to Problem $(MP)$ is fractional.

\subsection{Determinant of emission reduction potential}
\begin{figure}[htb!]
	\centering
  \includegraphics[width=0.85\textwidth]{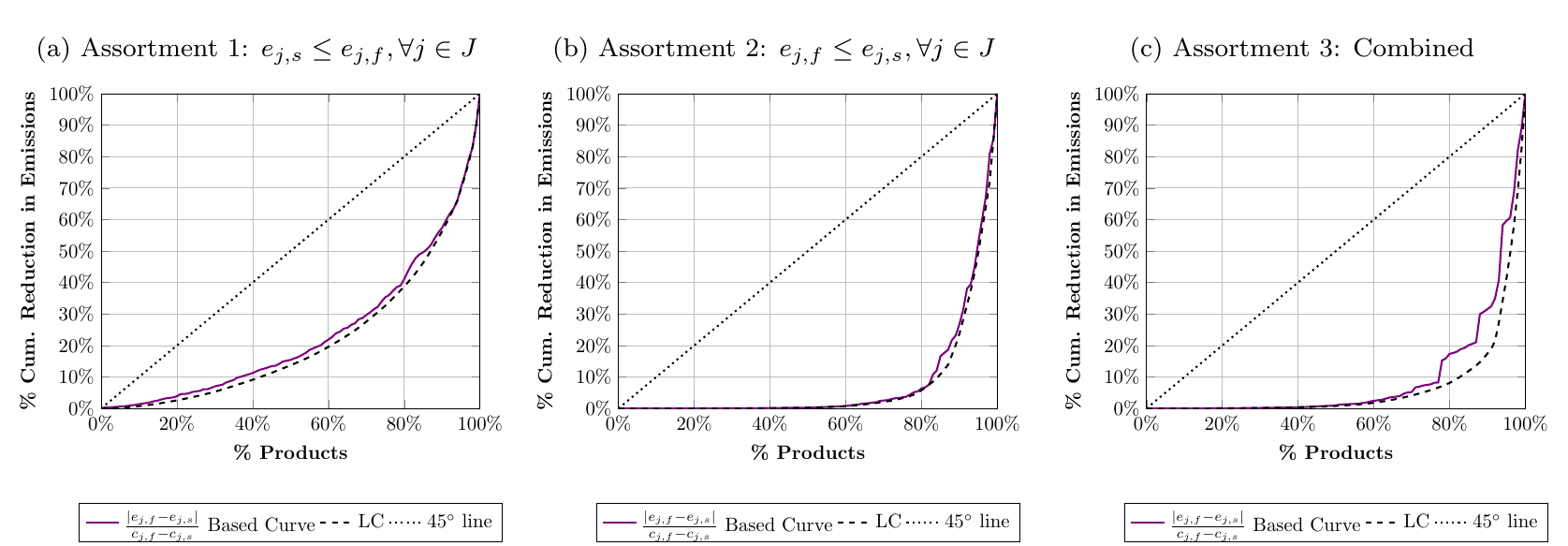}\\
	\caption{The cumulative emissions reduction share of items arranged based on different criteria.}
%%	Solid line - items are sorted in ascending order of the $\frac{|e_{j,f}- e_{j,s}|}{ c_{j,f}- c_{j,s}}$ ratio,\\
%	Dashed line - Lorenz curve for emissions reduction,\\
%	Dotted line - $45^{\circ}$ line representing the circumstances of an even distribution of emissions reduction.}}
  \label{fig:gini}
\end{figure}

The base case analysis of our DMS approach indicates that the emission reductions are not evenly distributed across the products. Products can be ordered by their contribution to emission reduction following the DMS optimization. In this manner we can construct the cumulative reduction in emissions as shown in the Lorenz curves (Figure \ref{fig:gini}) with the dashed line. 
Figure \ref{fig:gini} shows that 20\% of the products in assortments 1 through 3 account for 61.22\%, 94.19\%, and 91.88\% of the emission reduction, respectively.
This suggests that most of the emission reduction can be achieved by using dynamic mode selection for a limited subset of a given assortment. 
Although it is possible to determine the limited subset of products that account for most of the emission reduction after performing the DMS optimization, it would be convenient to know which products to focus on without having to solve a sophisticated optimization problem. Suppose we order products in increasing order of $\frac{|e_{j,f}- e_{j,s}|}{ c_{j,f}- c_{j,s}}$, i.e., we sort products according to the how much emission can be saved by using the least polluting transport mode relative to the additional cost of the faster transport mode. Figure \ref{fig:gini} shows the cumulative emission saving by products ordered this way with the solid line. Here we see that focusing on the 20\% of 
products in assortments 1 through 3 for which $\frac{|e_{j,f}- e_{j,s}|}{ c_{j,f}- c_{j,s}}$ is highest, already achieves 58.93\%, 93.63\%, and 82.64\% of the potential emission reduction, respectively. Thus firms seeking to minimize the carbon footprint of their inbound logistics should focus their attention on products for which the difference in emission in different transport modes is large relative to the additional cost of fast transportation modes. 
That is emission differences between modes relative to cost difference between modes is the main determinant of emission reduction potential for a given assortment. 

\subsection{Comparative statics}
\label{subsec:statics}
In this section, we study how changes in the input parameters with respect to the base case affect the performance of the blanket mode selection and the static mode selection approach. 
In what follows, we keep the emission target fixed at a 50 percent reduction of the reducible emissions, and we study the effects of changing a certain input parameter while generating the other input parameters as in the base case, i.e. as in Table \ref{tab:baseCase}. 
We also investigate the effect of scaling the emission differences between the most polluting and least polluting transport modes when the target on the total emissions is kept fixed. To achieve this, we first generate emission units as in the base case. We subsequently change the emission units of the most polluting transport mode through scaling $\mid e_{j,s} - e_{j,f} \mid$ by a constant $\delta_e$ while keeping the emission units from the least polluting mode fixed at its base level.   
%To do this, we generate emission units as in the base case, and subsequently increase the actual emission unit of the most polluting transport mode by scaling the difference with the least polluting transport mode by a constant $\delta_e$. Observe that $\delta_e$ equals $1$ in the base case.  
The changes in the parameters we investigate are summarized in Table \ref{tab:changes}.  

\begin{table}[htbp]
  \centering
  \caption{\textsf{Changes in the base case input parameters.}}
    \fontsize{8pt}{9pt}\selectfont
    \begin{tabular}{llll}
    \toprule
    \multicolumn{1}{l}{Parameter} &  Generation &  Base case  & \multicolumn{1}{l}{Changes}  \\
    \midrule 
    $\vert \mathnormal{J} \vert$ & & 100  & $\vert \mathnormal{J} \vert \in \{40, 60, 80\} $\\
    $\mu_{D_j}$  &    $\Gamma(100, CV_{\mu_{D_j}})$   & 0.5 & $CV_{\mu_{D_j}}\in\{ 0.3, 0.4,0.6,0.7\}$ \\
    $CV_{D_j}$  &    $B(0.9, 0.25, s_{D_j})$   & 0.3 & $s_{D_j} \in \{0.2, 0.25, 0.35, 0.4\}$ \\
 $\rho_{\mu_{D_j}, h_j}$  &       & -0.5 &  $\rho_{\mu_{D_j}, h_j}\in \{-0.3, -0.4, -0.6, -0.7\}$ \\
    $\psi_{p}$ &      & 9 &  $\psi_{p}  \in \{3, 4,5,19,99\}$ \\
    $\chi_j^c$ &  $B(0.25, \sigma_c)$    & 0.1 &  $\sigma_c  \in \{0.15,0.2,0.3,0.35\}$ \\
    $l_{j,s}$ &     & 3 &  $l_{j,s} \in \{2,4\}$ \\
    $\delta_e$ & & 1 & $\delta_e \in \{0.8,0.9,1.1,1.2\}$ \\ 
    \bottomrule
    \end{tabular}%
  \label{tab:changes}%
\end{table}%

\begin{figure}[htb!]
	\centering
  \includegraphics[width=0.8\textwidth]{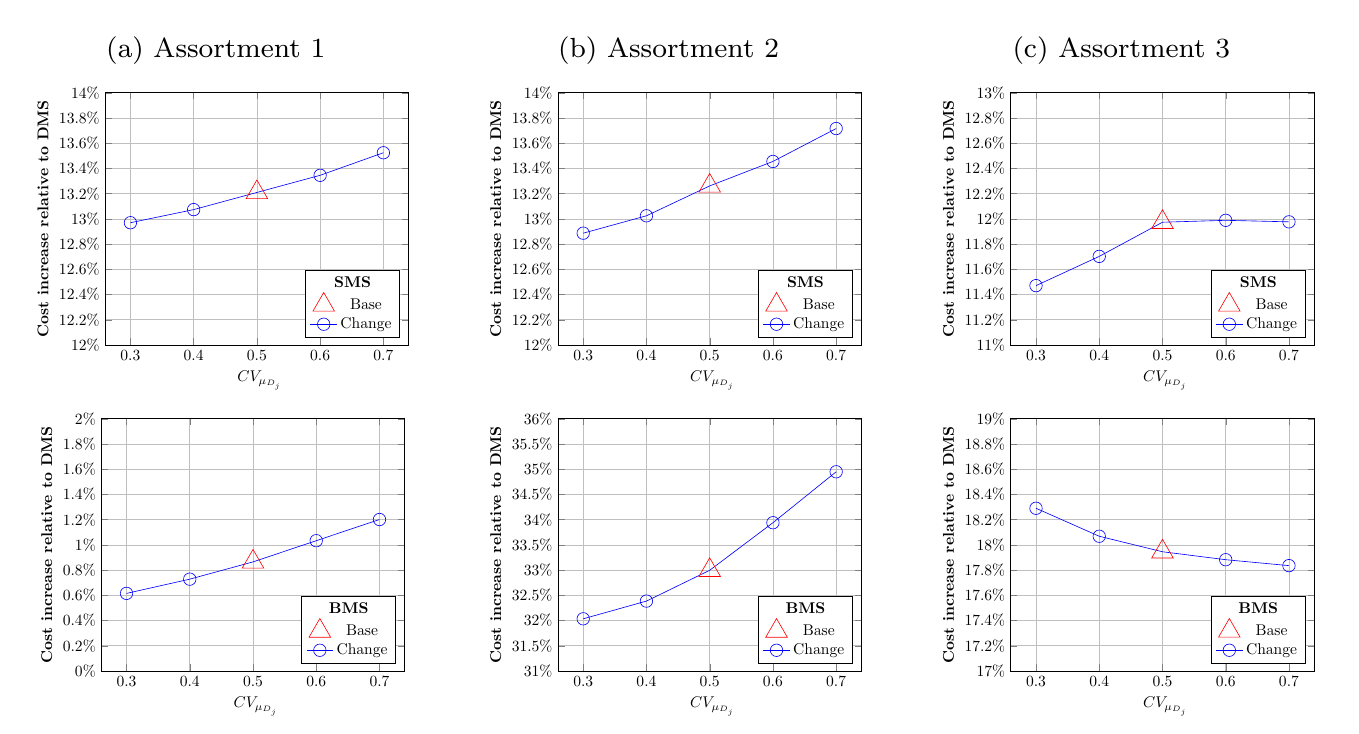}\\
	\caption{\textsf{Effect of changing $CV_{\mu_{D_j}}$ while keeping the rest of the parameters as in the base case.}}
  \label{fig:mu_d}
\end{figure}

\begin{figure}[htb!]
	\centering
  \includegraphics[width=0.8\textwidth]{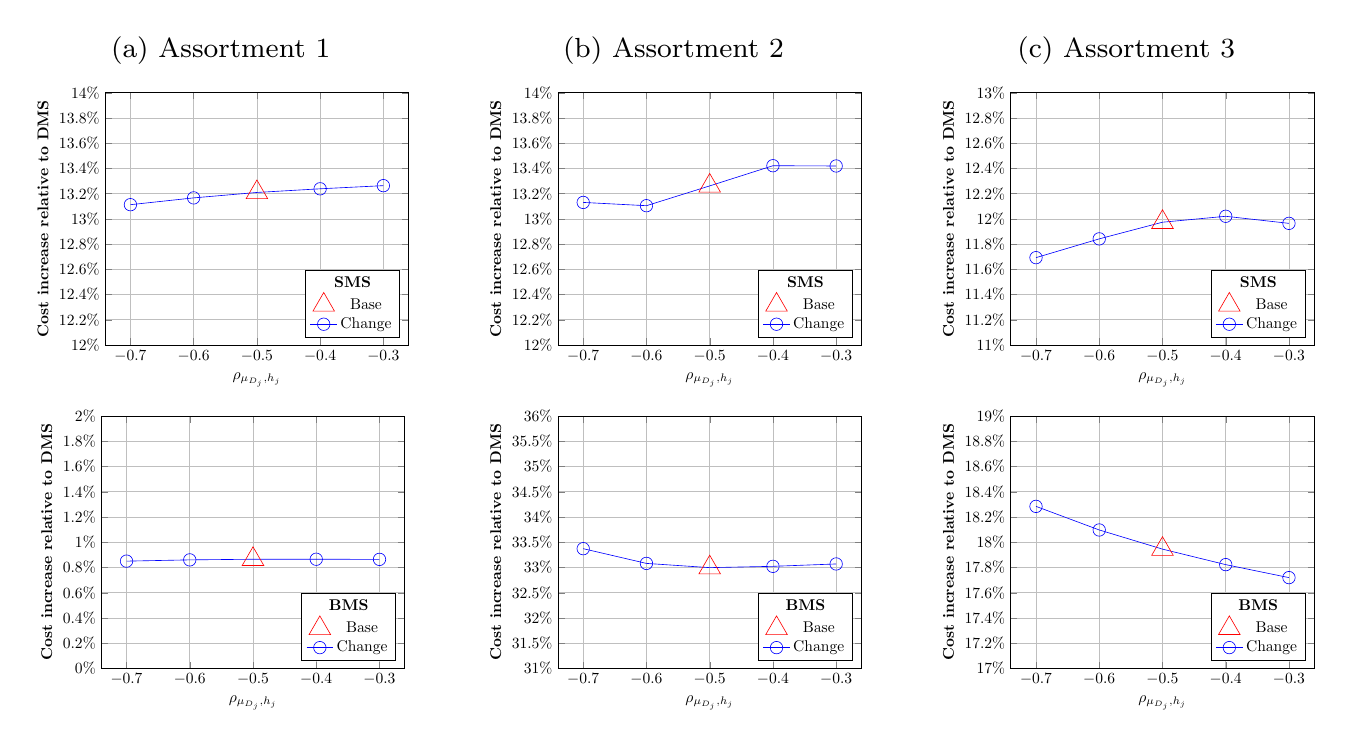}\\
	\caption{\textsf{Effect of changing $\rho_{\mu_{D_j}, h_j}$ while keeping the rest of the parameters as in the base case.}}
  \label{fig:rho}
\end{figure}

Figure \ref{fig:mu_d} shows the effect of changing the coefficient of variation of the gamma distribution from which we sample the mean demands per period for each product. 
The figure indicates that this effect is relatively limited. 
With respect to the base case, both $\%\mathit{BMS}$ and $\%\mathit{SMS}$ only change up to 1 percent point for all three assortment types.
We can draw a similar conclusion for the effect of changing the Pearson correlation coefficient between the holding cost and the mean demand per period of each product. 
Figure \ref{fig:rho} shows that both $\%\mathit{BMS}$ and $\%\mathit{SMS}$ change at most 1 percent point with respect to the base case for all three assortment types.

\begin{figure}[htb!]
	\centering
  \includegraphics[width=0.8\textwidth]{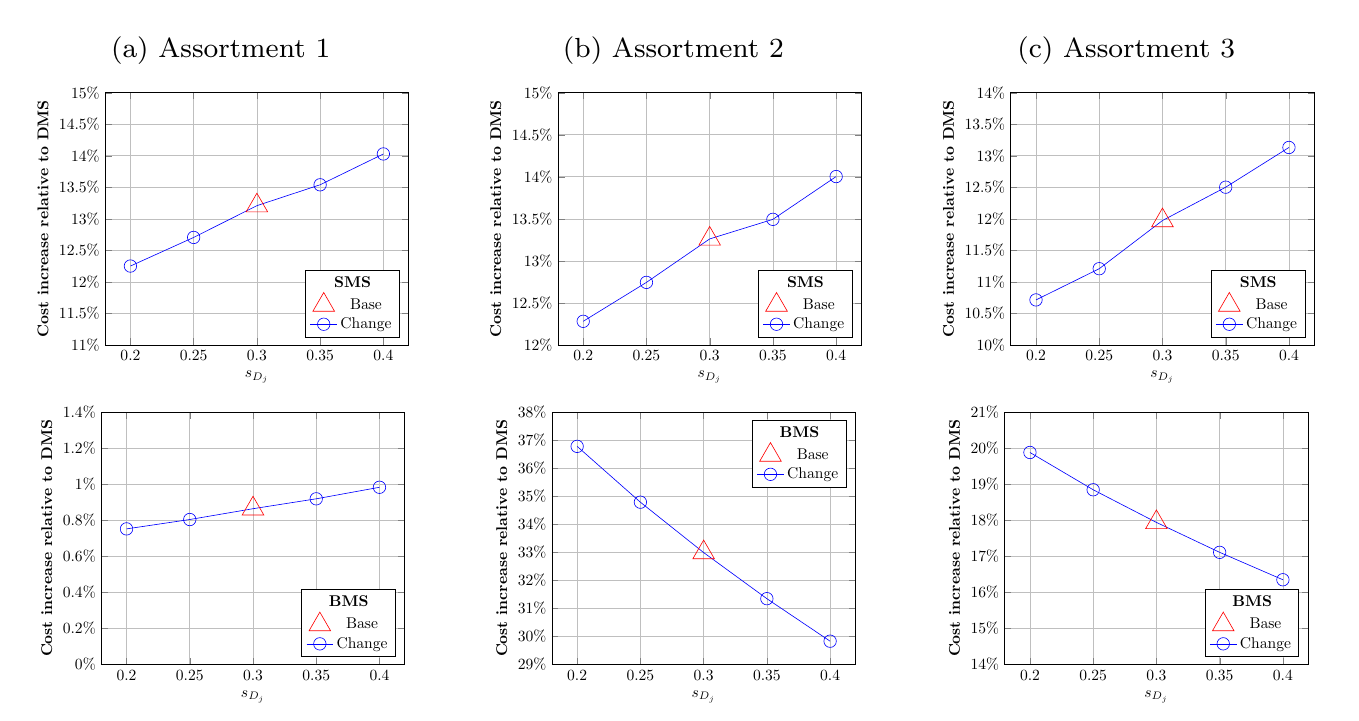}\\
	\caption{\textsf{Effect of changing $s_{D_j}$ while keeping the rest of the parameters as in the base case.}}
  \label{fig:s_d}
\end{figure}

Alterations in the shift parameter of the beta distribution from which we sample the coefficient of variation of the demand per period for each product has a relatively moderate effect on both $\%\mathit{BMS}$ and $\%\mathit{SMS}$. Figure \ref{fig:s_d} illustrates that for all assortment types, the $\%\mathit{SMS}$ tends to increase in the variability of the demand while the $\%\mathit{BMS}$ decreases.  
This indicates that the  flexibility to dynamically ship products with two transport modes has particular merit in highly variable demand settings.

\begin{figure}[htb!]
	\centering
  \includegraphics[width=0.8\textwidth]{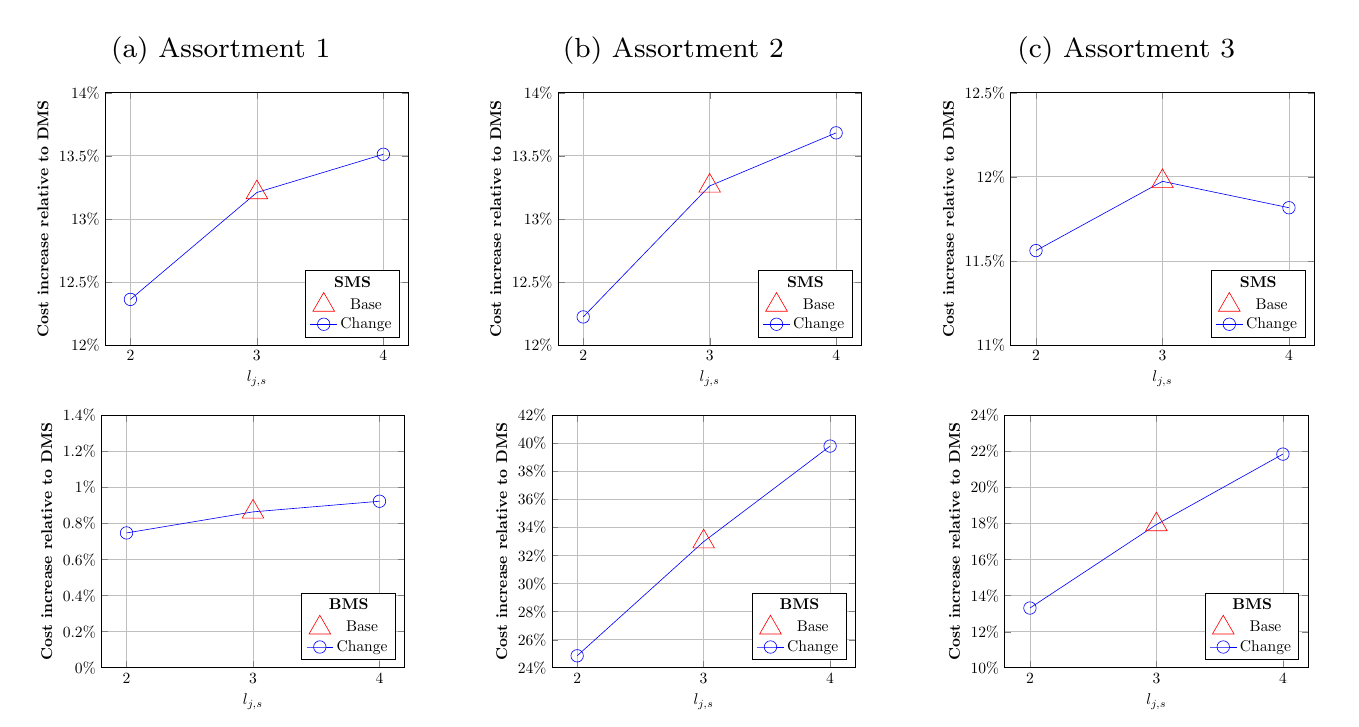}\\
	\caption{\textsf{Effect of changing $l_{j,s}$ while keeping the rest of the parameters as in the base case.}}
  \label{fig:l_s}
\end{figure}

Figure \ref{fig:l_s} indicates that for all assortment types, an increase (decrease) in the lead time difference between the fastest and the slowest transport modes of each product leads to an increase (decrease) in both $\%\mathit{SMS}$ and $\%\mathit{BMS}$. 
The blanket mode selection approach seems to be more susceptible to changes in the lead time difference than the static mode selection approach. For assortment type 2, for instance, an increase in the lead time difference to 4 leads to an increase in $\%\mathit{BMS}$ of 7 percent points with respect to the base case. By contrast, $\%\mathit{SMS}$ increases only slightly by 0.5 percent points. 
This can be attributed to the fact that the blanket mode selection approach imposes constraints on the emissions of each individual product while the static mode selection approach imposes a single constraint on the entire assortment of products. 

\begin{figure}[htb!]
	\centering
  \includegraphics[width=0.8\textwidth]{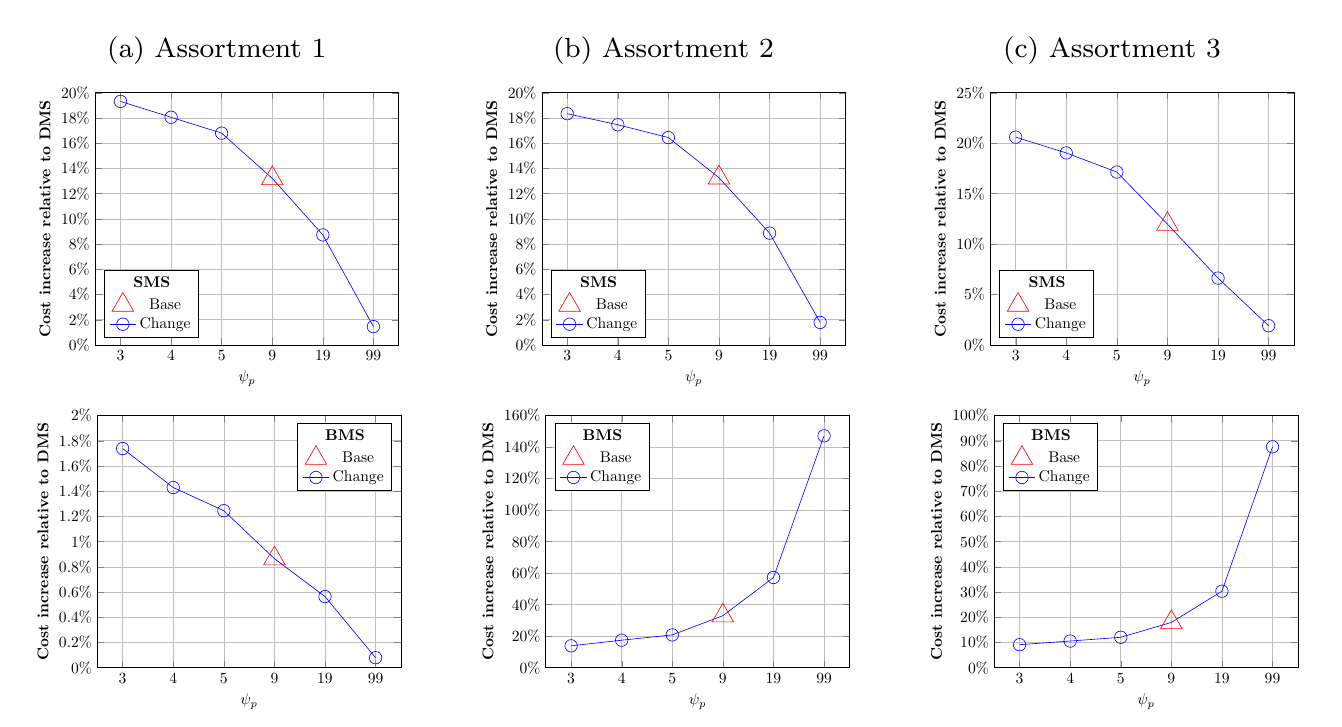}\\
	\caption{\textsf{Effect of changing $\psi_p$ while keeping the rest of the parameters as in the base case.}}
  \label{fig:psi_p}
\end{figure}

Figure \ref{fig:psi_p} illustrates the effect of changing the critical ratio for all products through varying $\psi_p$. We conclude that this effect is quite large.  For assortment type 1, for instance, $\%\mathit{SMS}$ varies from 20 percent to 2 percent. While $\%\mathit{SMS}$ seems to decrease in the critical ratio for all products, $\%\mathit{BMS}$ tends to increase. For assortment type 2 and 3, for instance, $\%\mathit{BMS}$ increases from 15 percent to over 80 percent. 
These effects can be explained by the fact that as the critical ratios of all products approach 1, our dynamic mode selection approach will mimic the static mode selection approach in which the assortment wide emission constraint is met by relying on less polluting transport mode for product for which this is relatively cheap to do so. The blanket mode selection approach, however, must meet emission targets for each product individually which leads to poor performance if we increase the critical ratios for all products. 

\begin{figure}[htb!]
	\centering
  \includegraphics[width=0.8\textwidth]{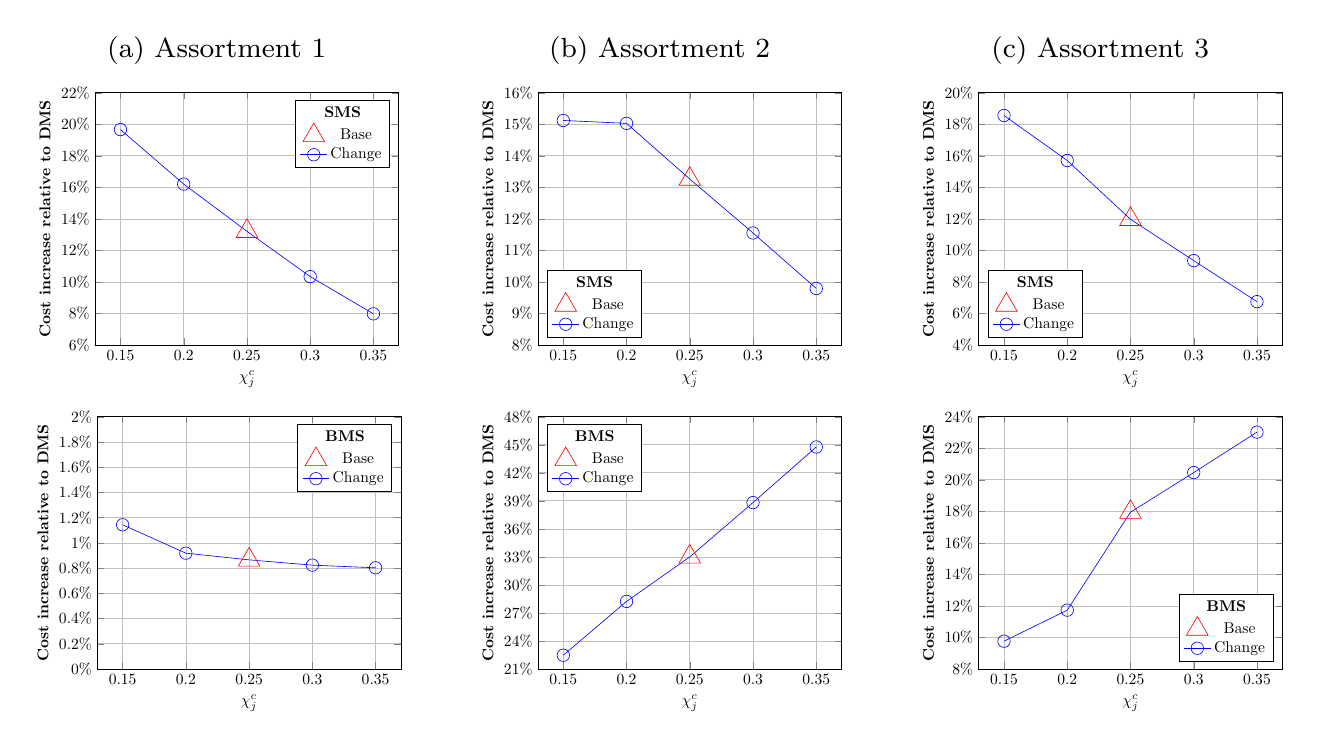}\\
	\caption{\textsf{Effect of changing $\chi_j^c$ while keeping the rest of the parameters as in the base case.}}
  \label{fig:chi_p}
\end{figure}

\begin{figure}[htb!]
	\centering
  \includegraphics[width=0.8\textwidth]{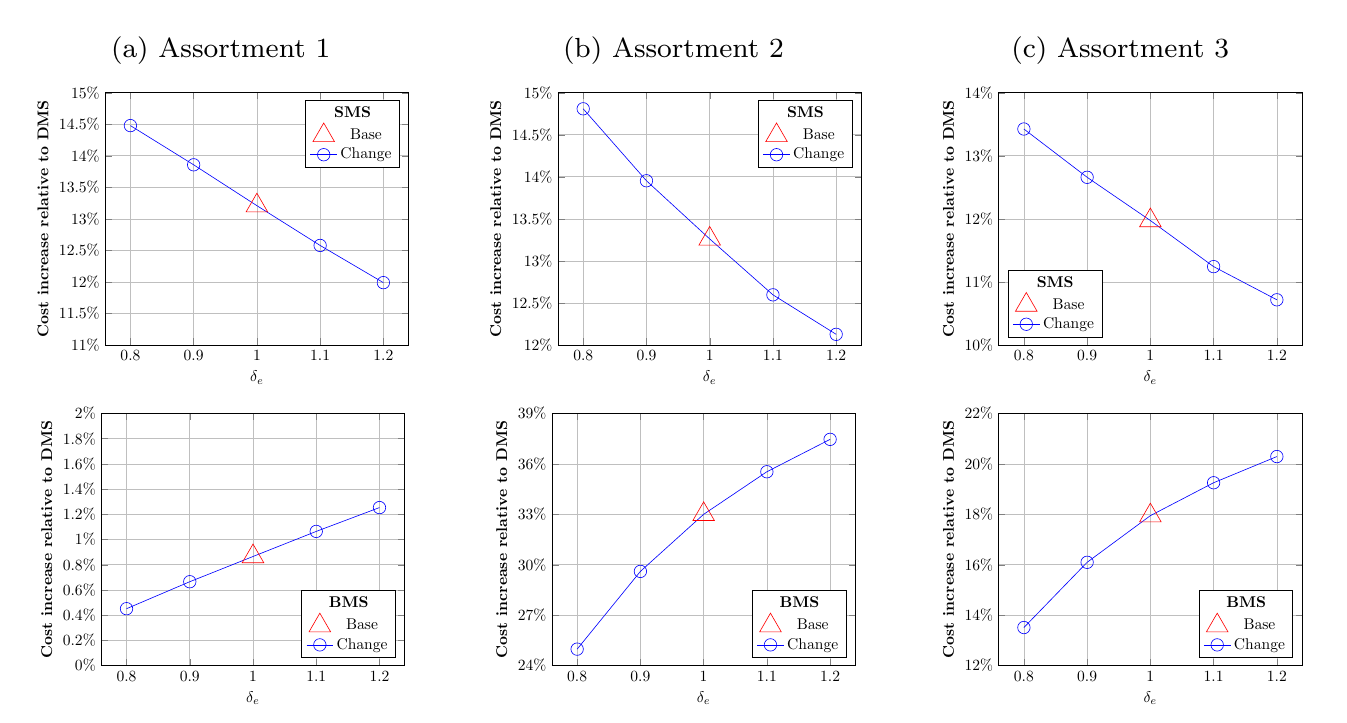}\\
	\caption{\textsf{Effect of changing $\delta_e$ while keeping the rest of the parameters as in the base case.}}
  \label{fig:delta_e}
\end{figure}

\begin{figure}[htb!]
	\centering
  \includegraphics[width=0.8\textwidth]{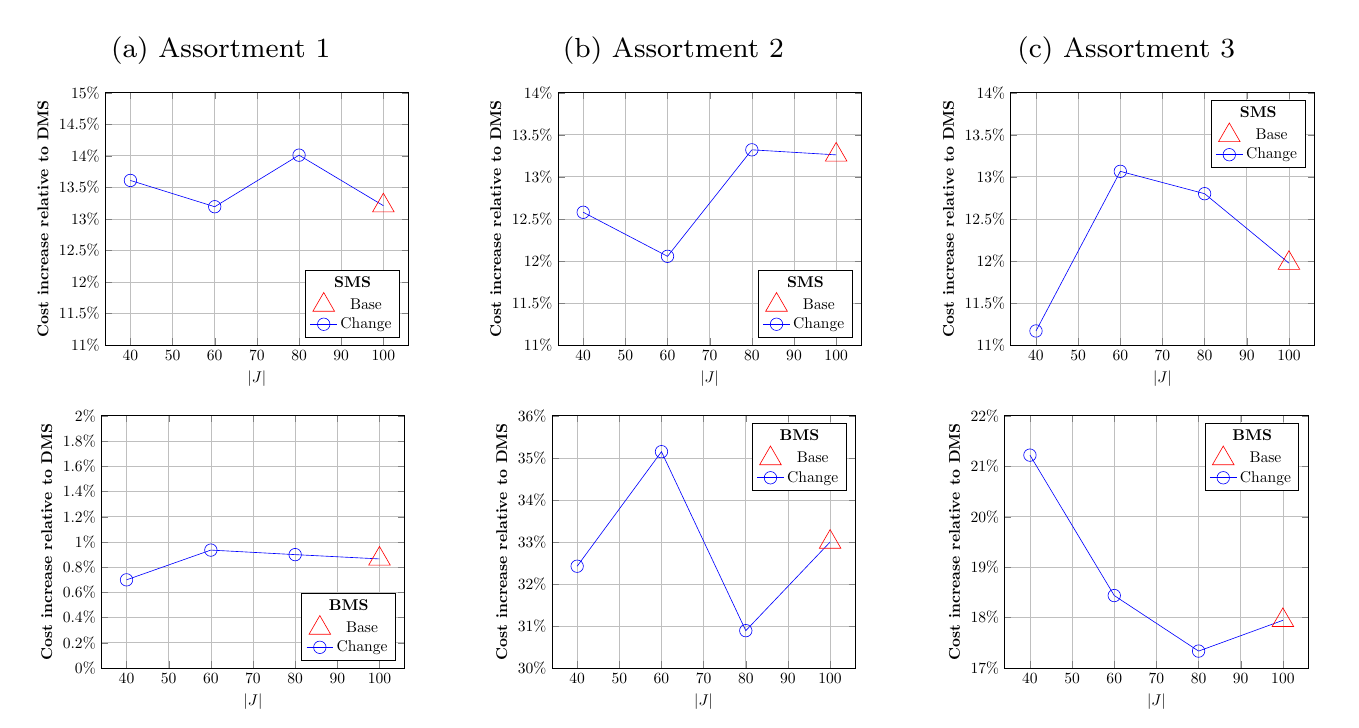}\\
	\caption{\textsf{Effect of changing $\vert J \vert$ while keeping the rest of the parameters as in the base case.}}
  \label{fig:size}
\end{figure}

Figure \ref{fig:chi_p} indicates that $\%\mathit{SMS}$ decreases  in  the  cost  of  the  fast  transport mode. This can again be explained by the fact that our dynamic mode selection approach will also rely more on the cheaper transport mode as the cost premium for the fast transport mode increases, and that consequently the gap with the static mode selection approach decreases. By contrast, $\%\mathit{BMS}$ increases in the cost of the fast transport mode for assortment type 2 and 3. This can be attributed to the fact that the blanket mode selection approach, contrary to the other two approaches, imposes itemized emission constraints and relying on the most expensive but least polluting transport mode is therefore inevitable. Note that this is not true for assortment type 1 because there the fast, expensive transport mode is also the most polluting transport mode. 
We can draw similar conclusions for the effects scaling the emission units of the most polluting transport mode, see Figure \ref{fig:delta_e}.   
Finally, Figure \ref{fig:size} indicates that the impact of the assortment size on both $\%\mathit{SMS}$ and $\%\mathit{BMS}$ is relatively limited. With respect to the base case, both $\%\mathit{BMS}$ and $\%\mathit{SMS}$ only change up to 2 percent point for all three assortment types. 

\section{Concluding remarks}
\label{sec:conclusionDM}
As carbon emissions from the transportation sector are projected to increase over the next decades, it is important for companies to rethink their supply chain strategies and explicitly incorporate carbon emissions into their decision making.
In this paper, we have studied the inbound transport and inventory management decision making for a company that sells an assortment of products. 
The company wishes to minimize inventory costs while keeping the total emissions from the inbound transport of the entire assortment below a certain target level.
Each product can be shipped using two distinct transport modes.
As each mode has its own merits, we have proposed a dynamic mode selection model that allows the company to ship products with either mode depending on when one mode is more favorable than the other. 
Since the optimal policy for dual transport mode problems are known to be complex, we have assumed that shipment shipment quantities for each product are governed by a dual-index policy. 
We have formulated the resulting decision problem as a mixed integer linear program that we have solved through a column generation solution procedure. This column generation procedure decomposes the complex multi-product problem into smaller sub-problems per product. These sub-problems are readily solved through a simple bisection search over Newsvendor type problems.

In an extensive computational experiment, we have compared the performance of our dynamic mode selection approach with two alternative approaches that are considered state of the art.
The first benchmark, static mode selection, lacks the flexibility to dynamically ship products with two transport modes; it rather selects one transport mode for each product a priori.
The second benchmark, blanket mode selection, does have the flexibility to rely on two transport modes simultaneously but it makes transport decisions for each product individually rather than holistically for the entire assortment. 
%Our computational experiment considered three different types of assortments of products; one in which the slowest transport modes of all products are also the least polluting modes, one in which the fastest transport modes of all products are also the least polluting modes, and one in which for any product either the fastest or slowest transport mode is the least polluting mode.
Our computational experiments indicate that the value of our dynamic mode selection approach over the blanket mode selection approach is particularly high for assortments of products for which the fastest transport modes are not necessarily the most polluting transport modes. For such settings, our dynamic mode selection approach can reduce the long run average costs by 40 percent under moderate carbon emission targets. 
These huge savings can be attributed to the portfolio effect inherent to our approach. The computational experiments further indicate that dynamic mode selection can significantly outperform static mode selection. Under moderate emission targets, dynamically relying on two transport modes rather than a single transport mode can lead to cost savings of up to 15 percent.

Future studies can extend the current model by studying other settings with multiple transport modes such as multiple echelons in a serial system \citep[e.g.][]{lawson2000,ARTS2013} or assembly systems \citep[e.g.][]{angelus2016}. Alternatively one may consider more sophisticated dual mode heuristic policies such as the projected expedited inventory position policy \citep{MelvinArts}, capped base-stock policy \citep{sun2019robust}, or vector base-stock policy \citep{sheopuri2010new}.

% Acknowledgments here
\ACKNOWLEDGMENT{
The research of the first author is supported by the National Research Fund of Luxembourg through AFR grant 12451704. The authors gratefully acknowledge Ranit Sinha for his comments on earlier versions of the paper.  
}

\bibliographystyle{informs2014}
\bibliography{bib}

\begin{APPENDICES}

%********Finish**********

\section{Carbon accounting}
\label{sec:carbon}
In this section, we briefly explain how we determine the distribution functions that we use for pseudo-random generation of the unit emissions in our computational experiments. We utilize the \cite{UNComtrade} to calculate the average unit weights for 122 groups of products imported by The Netherlands in 2020. These product groups consist of two categories: (i) apparel goods that are imported from Vietnam and (ii) industrial goods that are imported from China and Germany. We consider air transport (from Tan Son Nhat international airport) and sea transport (from Haiphong port) as the fast and slow transport mode for the apparel category, respectively. For the industrial goods category, we assume sea transport from Shanghai, China, as the slow mode and road transportation from Stuttgart, Germany, as the fast mode.  

We rely on the Network for Transport Measures methodology \citep{NTMframework} to model and measure transportation emissions based on equation \ref{eq:emission}. This model has been widely used in literature \citep[e.g.,][]{Hoen,Hoen2014Switching}. Following the NTM methodology, we first compute the overall carbon emissions generated by a single vehicle and then allocate a proportion of those emissions to each freight unit carried by the vehicle. 

\textbf{Sea transportation.} We assume that all products are shipped via container. The average age of the container fleet worldwide is around 12 years and the average vessel size (dwt) of container ships with age 10-14 is 43,993 ton \citep{UNCTADreport}. Based on section 7 of the NTM framework and resolutions of the Marine Environment Protection Committee \citep{MEPC20362}, we approximate sea transportation emissions in  kilograms of $CO_2$ of one unit of a certain product with weight $w$ (in kilogram) for a certain trip with distance $d$ (in kilometers) using the following relation,
\[
    e_{sea}=w \cdot EI_{ship} \cdot 10^{-3} \cdot d
\]
where $EI_{ship}$ is kilograms of $CO_2$ emissions per kilogram weight per kilometer. Furthermore $EI_{ship}$ is computed through,
\[
    EI_{ship} = \frac{(a \cdot dwt^{-c} ) / (PDR_{ship} \cdot LCU)}{1.852}
\]
where $a$ and $c$ are constants, $dwt$ is the deadweight tonnage of the ship, $LCU$ is average load capacity utilization, $F (LCU)$ is fuel consumption as a function of load, and $PDR_{ship}$ is the payload of the ship. $1.852$ is the nautical mile to km conversion coefficient. For a container ship, NTM methodology states: $a=0.17422, c= 0.201, LCU=0.70, F(LCU)=1$, and $PDR=0.8$. Succinctly, we have for the total emissions in kilograms of $CO_2$ of one unit of a product with weight $w$ (in kilogram) for a sea trip with distance $d$ (in kilometers)
\begin{equation} \label{eq:seaCo}
    e_{sea} = w(1.996 \cdot 10^{-5} \cdot d).
\end{equation}
\textbf{Air transportation.} Our calculations for the emissions of air transportation are based on section 8 of the NTM Framework. We consider an Airbus A310-300 F as the aircraft. Based on the May 2021 Air Cargo Market Analysis of The International Air Transport Association, we assume an average international cargo load factor of 65\%. Following the NTM Framework, we have the following relation for air transportation emissions
\[
    e_{air} = \frac{w}{c_{max}}(CEF + VEF \cdot d),
\]
where $c_{max}$ is the maximum freight load, $CEF$ is the constant emissions factor, and $VEF$ is the variable emissions factor. $CEF$ and $VEF$ are the outcomes of applying a linear regression on real data provided by the NTM. We obtain the $CEF$ and $VEF$ parameters via interpolation over the associated tables provided by the NTM. We furthermore assume $c_{max}=39,000 kg$ as per section 8.3.1 and perform the interpolation on table 4.1 of section 8.2.1. Succinctly, we have for the total emissions in  kilograms of $CO_2$ of one unit of a product with weight $w$ (in kilogram) for an air trip with distance $d$ (in kilometers)
\begin{equation} \label{eq:airCo}
    e_{air} =w (1.525 \cdot 10^{-1} + 4.938 \cdot 10^{-4} \cdot d).
\end{equation}

\textbf{Road transportation.} We rely on \cite{Hoen} to obtain the emission units of road transportation. They too rely on the NTM framework to estimate $CO_2$ emissions from road transportation in Europe. In particular, they approximate the total emissions in  kilograms of $CO_2$ of one unit of a product with weight $w$ (in kilogram) for a truck trip with distance $d$ (in kilometers) as
\begin{equation} \label{eq:roadCo}
    e_{road} = w(3.214 \cdot 10^{-4} + 4.836 \cdot 10^{-5} \cdot d).
\end{equation}

\textbf{Distances.} We use NTMCalc Basic 4.0 \citep{NTMCalc}, which is an online tool provided by the NTM for approximating emissions, to calculate travel distances between the origin destination pairs as described at the beginning of this section. Based on this tool we find that the sea distance between Haiphong and Rotterdam is 9,610 nautical miles (17,798 km) and the sea distance between Shanghai and Rotterdam is 10,525 nautical miles (19,492 km). The distance traveled by aircraft between Tan Son Nhat international airport and Rotterdam The Hague Airport is 10,073 km, and the road distance between Stuttgart and Rotterdam is 633 km. With these distances, we compute the total kilogram $CO_2$ emissions of one unit of product with weight $w$ for each mode-trip category using Equations \eqref{eq:seaCo}-\eqref{eq:roadCo}. We call these emission coefficients and they are presented in Table \ref{tab:coefficients} below.

\begin{table}[ht]
\centering
\caption{\textsf{Emissions coefficients for each mode-trip category.}}
\label{tab:coefficients}
	\fontsize{8pt}{9pt}\selectfont
\begin{tabular}{l l l}
\toprule
Industry    & Slow Mode             & Fast Mode
\\ \midrule
Apparel    & $3.552 \cdot 10^{-1}$ & $5.127$\\
Industrial & $3.891 \cdot 10^{-1}$ & $3.093 \cdot 10^{-2}$
\\ \toprule
\end{tabular}
\end{table}

\textbf{Fitting distribution functions.} We use the emission coefficients from Table \ref{tab:coefficients} to calculate the unit emissions for the 122 groups of products mentioned at the beginning of this section. We subsequently use maximum likelihood estimation on the resulting emission units to find distribution functions from which we can sample the emission units of the fast and slow transport modes for all three assortment types in our computational experiment. The emission unit distributions for assortment type 1 are based on the apparel category, the emission unit distributions for assortment type 2 are based on the industrial category, and the emission unit distributions for assortment type 3 are based on both categories. The final distribution functions are provided in Table \ref{tab:emissionUnits}.

\section{Generating correlated random numbers}
\label{sec: random_numbers}
Suppose $X$ and $Y$ are two real random variables with marginal distribution functions $F$ and $G$, respectively. Suppose their joint distribution is bi-variate standard normal $\mathcal{N}_{\rho}$ with Pearson's correlation coefficient $\rho=\cov(X,Y)/(\sqrt{\Var[X]}\sqrt{\Var[Y]})$ \citep{Nelsen2006}. Let $Z$ be a vector of size two with independent random elements that have standard normal distributions $\Phi$, and let $W=AZ$ be a linear combination of $Z$ with \[
A=
    \begin{bmatrix}
    1 & 0\\
    \rho & \sqrt{1-\rho ^ 2}
    \end{bmatrix}.
\] 

It can be shown that $W$ has a bivariate normal distribution $\mathcal{N}_{\rho}$ with covariance matrix $\Sigma = AA^T$ \citep[see, e.g.,][]{Gut2009}. We use this result to sample from $X$ and $Y$ as follows:
\begin{enumerate}
    \item Generate the vector $Z=
            \begin{bmatrix}
                Z_1\\
                Z_2
            \end{bmatrix}$ by independently sampling from a standard normal distribution function,
    \item Calculate the bivariate normal sample $W=
            \begin{bmatrix}
                W_1\\
                W_2
            \end{bmatrix}=AZ$,
    \item Generate the required samples by inversion $
            \begin{bmatrix}
                X\\
                Y
            \end{bmatrix} =
            \begin{bmatrix}
                F^{-1}(\Phi(W_1))\\
                G^{-1}(\Phi(W_2))
            \end{bmatrix}$.
\end{enumerate}

\section{Benchmark approaches}
\label{sec:BM_approaches}
 The mathematical formulation for the blanket mode selection approach, which enforces emission constraints $\mathcal{E}_j^{max}$ per item $j\in J$, is called Problem $(\mathit{BMS})$ and is given as follows:
 \begin{alignat*}{5}
&(\mathit{BMS})\qquad \qquad  	  &&\displaystyle \min       &&\quad && C(\mathbf{S}_f,\mathbf{\Delta}) &&  \\
&        				      &&\text{subject to} &&      \quad &&\mathnormal{E}_j(S_{j,f},\Delta_j) \leq \mathcal{E}_j^{max}, &&\qquad \forall j\in J, \\
&											&&									&&	\quad		&&(S_{j,f}, \Delta_j) \in (\mathbb{R} \times \mathbb{R}_0), \quad  &&\qquad \forall j \in J. 
 \end{alignat*}

 The mathematical formulation for the static mode selection approach, which selects one transportation mode per item, is called Problem $(\mathit{SMS})$ and is given as follows:
\begin{alignat*}{5}
&(\mathit{SMS})\qquad\qquad && \min && \quad && \sum_{j\in J} C_{j,s}(S^*_{j,s}) x_{j,s} + \sum_{j\in J} C_{j,f}(S^*_{j,f})x_{j,f}  \\
& && \text{subject to} && \quad && \sum_{j\in J} E_{j,s}x_{j,s} +\sum_{j\in J} E_{j,f}x_{j,f} \leq \mathcal{E}^{max}, && \\
& && &&  &&  x_{j,s},x_{j,f}\in\{0,1\},  && \forall j \in J,
\end{alignat*}
where $C_{j,s}(S_{j,s}):=h_j \E\left[\left(S_{j,s}-\textstyle\sum_{t=0}^{l_{j,s}}D^t_j\right)^+\right]+ p_j \E\left[\left(\textstyle\sum_{t=0}^{l_{j,s}} D^t_j - S_{j,s}\right)^+\right]+c_{j,s}\E[D_j]$, and
$C_{j,f}(S_{j,f}):=h_j \E\left[\left(S_{j,f}-\textstyle\sum_{t=0}^{l_{j,f}}D^t_j\right)^+\right]+ p_j \E\left[\left(\textstyle\sum_{t=0}^{l_{j,f}} D^t_j - S_{j,f}\right)^+\right]+c_{j,f}\E[D_j]$
are the average cost-rate for using the slow and fast transport mode respectively,
$S^*_{j,s}:=\argmin_{S_{j,s}}C_{j,s}(S_{j,s})$, and $S^*_{j,f}:=\argmin_{S_{j,f}}C_{j,f}(S_{j,f})$ are the base-stock levels that minimize inventory related costs when exclusively using the slow and fast transport mode respectively, and $E_{j,s}:=e_{j,s}\E[D_j]$ and $E_{j,f}:=e_{j,f}\E[D_j]$ are the emissions per time unit of shipping exclusively with the slow and fast transport mode respectively.

\end{APPENDICES}

\end{document}